\documentclass[12pt,reqno]{NumPDEsArticle}

\usepackage[utf8]{inputenc}
\usepackage[T1]{fontenc}
\usepackage[english]{babel}
\usepackage{fullpage}
\usepackage{enumitem}
\usepackage{amsmath,amssymb}
\usepackage[dvipsnames]{xcolor}
\usepackage{array}
\usepackage{multirow}
\usepackage{mathtools}
\usepackage{xspace}
\usepackage{todonotes}
\usepackage{graphicx}
\usepackage{hyperref}

\usepackage{tikz}
\usepackage{pgfplots}
\usepackage{pgfplotstable}
\usetikzlibrary{calc,angles,quotes}
\usetikzlibrary{arrows.meta,positioning}
\usetikzlibrary{matrix}
\usepackage{pgf-umlcd}

\hypersetup{
	colorlinks,
	citecolor=OliveGreen,
	linkcolor=BrickRed,
	urlcolor=NavyBlue
}

\usepackage[backend=biber,style=trad-alpha,safeinputenc]{biblatex}
	
\AtEveryBibitem{%
	\clearfield{issn}%
	\clearfield{isbn}%
}

\DeclareFieldFormat{titlecase}{#1}
	
\usepackage{listings}

\definecolor{codegray}{rgb}{0.5,0.5,0.5}
\definecolor{backcolour}{rgb}{0.95,0.95,0.92}
\definecolor{matlabblue}{RGB}{14,0,255}
\definecolor{matlabgreen}{RGB}{2,128,9}
\definecolor{matlabpurple}{RGB}{170,4,249}

\definecolor{color1}{HTML}{0072BD}
\definecolor{color2}{HTML}{D95319}
\definecolor{color3}{HTML}{EDB120}
\definecolor{color4}{HTML}{7E2F8E}
\definecolor{color5}{HTML}{77AC30}

\makeatletter
\lstdefinestyle{mystyle}{
	language=MATLAB,
	commentstyle=\color{matlabgreen},
	keywordstyle=\color{matlabblue},
	numberstyle=\tiny\color{codegray},
	stringstyle=\color{matlabpurple},
	basicstyle=\normalfont\ttfamily\lst@ifdisplaystyle\scriptsize\fi,
	breakatwhitespace=false,
	breaklines=true,
	captionpos=b,
	keepspaces=true,
	numbers=left,
	numbersep=5pt,
	showspaces=false,
	showstringspaces=false,
	showtabs=false,
	xleftmargin=\parindent,
	tabsize=2
}
\makeatother

\defBlackboardCharacter{N}
\defBlackboardCharacter{R}
\defBlackboardCharacter{T}
\defBlackboardCharacter{Z}

\defCaligraphicCharacter{C}
\defCaligraphicCharacter{E}
\defCaligraphicCharacter{I}
\defCaligraphicCharacter{N}
\defCaligraphicCharacter{M}
\defCaligraphicCharacter{P}
\defCaligraphicCharacter{Q}
\defCaligraphicCharacter{S}
\defCaligraphicCharacter{T}
\defCaligraphicCharacter{U}
\defCaligraphicCharacter{V}
\defCaligraphicCharacter{X}

\newcommand{\nC}{n_{\VV}}
\newcommand{\nS}{n_{\EE}}
\newcommand{\nE}{n_{\TT}}

\newcommand{\matlab}{\textsc{Matlab}\xspace}
\newcommand{\moofem}{MooAFEM\xspace}
\newcommand{\cpp}{C\kern-.05em\raisebox{.3ex}{\textbf{\tiny++}}\xspace}

\newcommand{\coarse}{H}
\newcommand{\fine}{h}
\newcommand{\refine}{{\tt refine}}
\newcommand{\mat}[1]{\mathbf{#1}}
\newcommand{\A}{\mat{A}}

\renewcommand{\vec}[1]{\boldsymbol{#1}}
\newcommand{\normalvec}{\vec{n}}

\renewcommand{\d}[1]{\,\mathrm{d}#1}
\let\div\relax
\DeclareMathOperator{\div}{div}
\DeclareMathOperator{\conv}{conv}
\newcommand{\D}{\mathrm{D}}

\newcommand{\norm}[2][]{#1\| #2 #1\|}
\newcommand{\set}[3][]{#1\{ #2 \, #1| \, #3 #1\}}

\newcommand{\enorm}[2][]{#1|\!#1|\!#1| #2 #1|\!#1|\!#1|}
\newcommand{\jump}[1]{[\![#1]\!]}

\let\additionalLink\footnote

\title{MooAFEM: An object oriented Matlab code for higher-order adaptive FEM for (nonlinear) elliptic PDEs}

\author{Michael Innerberger}
\address{TU Wien, Institute of Analysis and Scientific Computing, Wiedner Hauptstr.~8--10/E101/4, 1040 Wien, Austria}
\email{Michael.Innerberger@asc.tuwien.ac.at \qquad\rm (corresponding author)}
\email{Dirk.Praetorius@asc.tuwien.ac.at}

\author{Dirk Praetorius}

\thanks{The authors acknowledge support through the Austrian Science Fund (FWF) through the doctoral school \emph{Dissipation and dispersion in nonlinear PDEs} (grant W1245), the special research program \emph{Taming complexity in PDE systems} (grant SFB F65), and the standalone project \emph{Computational nonlinear PDEs} (grant P33216).}

\date{\today}

\begin{document}
	
\begin{abstract}
	We present an easily accessible, object oriented code (written exclusively in \matlab) for adaptive finite element simulations in 2D.
	It features various refinement routines for triangular meshes as well as fully vectorized FEM ansatz spaces of arbitrary polynomial order and allows for problems with very general coefficients.
	In particular, our code can handle problems typically arising from iterative linearization methods used to solve nonlinear PDEs.
	Due to the object oriented programming paradigm, the code can be used easily and is readily extensible.
	We explain the basic principles of our code and give numerical experiments that underline its flexibility as well as its efficiency.
\end{abstract}

\maketitle


\section{Introduction}\label{sec:intro}

In today's ever-changing landscape of mathematical software, \matlab has successfully reinforced its role as a \textsl{de facto} standard for the development and prototyping of various kinds of algorithms for numerical simulations.
In particular, it has proven to be an excellent tool for academic education, e.g., in the field of the partial differential equations (PDEs), which can be solved by the finite element method (FEM); see, e.g., \cite{lichen-comppde}.
While the impact of vectorization for efficient FEM codes in MATLAB is well-known in the literature \cite{ifem,p1afem,csv2019,mv2022,varfem}, we show, on the one hand, how FEM codes in MATLAB can massively benefit from object oriented programming (OOP) and, on the other hand, how OOP even simplifies the efficient implementation of FEM for nonlinear PDEs.

In this work, we present our own freely available code \moofem (\matlab object oriented adaptive FEM)~\cite{code}.
It is specifically tailored to (higher-order) adaptive FEM simulations involving linear second-order elliptic PDEs with very general coefficients:
Let $\Omega \subset \R^2$ be a bounded domain with polygonal boundary $\partial \Omega$ that is split into \emph{Robin}, \emph{Neumann}, and \emph{Dirichlet} boundary, i.e., $\partial \Omega = \overline{\Gamma}_R \cup \overline{\Gamma}_N \cup \overline{\Gamma}_D$ and $\Gamma_R$, $\Gamma_N$, $\Gamma_D$ are pairwise disjoint.
The model problem then reads 
\begin{subequations}
\label{eq:model}
\begin{align}
	\label{eq:model-pde}
	-\div \A \nabla u + \vec{b} \cdot \nabla u + c u
	&=
	f - \div \vec{f}
	&& \hspace{-.1\textwidth} \text{in } \Omega,\\
	\label{eq:model-robin}
	\alpha \, u + \A \nabla u \cdot \normalvec
	&=
	\gamma
	&& \hspace{-.1\textwidth} \text{on } \Gamma_R,\\
	\label{eq:model-neumann}
	\A \nabla u \cdot \normalvec
	&=
	\phi
	&& \hspace{-.1\textwidth} \text{on } \Gamma_N,\\
	\label{eq:model-dirichlet}
	u
	&=
	0
	&& \hspace{-.1\textwidth} \text{on } \Gamma_D.	
\end{align}
\end{subequations}
Note that also inhomogeneous Dirichlet data $g \in H^{1/2}(\Gamma_D)$ can be handled by the usual superposition ansatz, i.e., $u = u_D + \widehat{g}$, where $\widehat{g} \in H^1(\Omega)$ is a lifting of the inhomogeneous boundary data (e.g., by $L^2(\Omega)$-projection or by nodal interpolation in the discrete case) and $u_D \in H^1(\Omega)$ with $u_D = 0$ on $\Gamma_D$ solves an appropriate variant of~\eqref{eq:model}; see, e.g.,~\cite{ern-guermond,bcd2004,axioms,fpp2014}.
Due to the modular design of \moofem, also other methods for incorporating inhomogeneous Dirichlet data (such as penalty methods~\cite{babuska1973} or Nitsche's method~\cite{nitsche1971}) can be implemented analogously to Robin-type boundary conditions.

\moofem is able to discretize problem~\eqref{eq:model} with conforming FEM spaces of arbitrary polynomial order and, by use of OOP, allows that the coefficients $\A$, $\vec{b}$, $c$, and $\alpha$, as well as the data $f$, $\vec{f}$, $\gamma$, and $\phi$ can be any function that depends on a spatial variable; in particular, FEM functions are also valid as coefficients.
In contrast to existing object oriented \matlab FEM codes like \cite{ofem}, this allows to treat also non-linear PDEs (as well as nonlinear boundary conditions), since these must be solved by iterative linearization techniques in practice, which lead to problems of the form~\eqref{eq:model}, where the coefficients as well as the data depend on previous iterates (see Section~\ref{subsec:nonlinear} below).
The implementation of such (nonlinear) problems with our code is fairly easy and extensively documented below.
Furthermore, all computations (e.g., function access, quadrature, and assembly) are vectorized efficiently such that the measured computation time scales optimally with the number of degrees of freedom.

Through proper referencing mechanics enabled by the use of OOP in \matlab, all vectorizations are hidden behind an application layer.
Keeping class design lean and hierarchies flat, together with consistent, descriptive naming and an extensive documentation, this application layer is highly readable while still granting the user full control over the underlying algorithms (e.g., quadrature rules and interpolation methods).
Furthermore, separation into several loosely coupled modules, an extensive unit-test suite covering the critical parts of the code, and well-designed interfaces for mathematical operations allow for facile extension of our code.
Thus, \moofem combines the flexibility and reusability of object oriented software with the accessibility of the high-level scripting language \matlab, providing a much needed comprehensive, flexible, and efficient library for AFEM research with only about one thousand executable lines.

\subsection*{Software requirements}
The \moofem package requires an installation of \matlab version R2020b (v9.9) or later.
We stress that, besides the core \matlab libraries, no other toolboxes are required.
For this reason, it may also be ported to Octave with reasonable effort (although this was never attempted).

\subsection*{Outline}

We first comment on the precise scope of our software package and reinforce the arguments in favor of object oriented design in Section~\ref{sec:justification}.
In particular, we give a schematic algorithm for AFEM that showcases the natural syntax of the \moofem package.
The general structure of the code is outlined in Section~\ref{sec:code-structure}: an exhaustive overview over available features is given and the most important ones are explained in more depth.
In Section~\ref{sec:data-structures}, we explain the underlying memory layout of the data structures used in \moofem and other implementational aspects.
Finally, we extend the basic AFEM algorithm from Section~\ref{sec:justification} to more interesting examples in Section~\ref{sec:examples}: higher-order FEM, goal-oriented FEM, and iterative linearization of nonlinear equations.


\section{Adaptive algorithm and importance of OOP}\label{sec:justification}

\subsection{Adaptive FEM}

To solve equation~\eqref{eq:model}, we employ FEM with underlying triangulation $\TT_\coarse$ of $\Omega$.
To this triangulation, we associate the FEM space $\SS^p(\TT_\coarse) := \PP^p(\TT_\coarse) \cap H^1(\Omega)$ with
\begin{equation*}
	\PP^p(\TT_\coarse)
	:=
	\set[\big]{v \in L^2(\Omega)}{v|_T \text{ is a polynomial of degree $p$ for all } T \in \TT_\coarse}.
\end{equation*}
With $H^1_D(\Omega) := \set{v \in H^1(\Omega)}{v|_{\Gamma_D} = 0}$, we set $\SS^p_D(\TT_\coarse) := \SS^p(\TT_\coarse) \cap H^1_D(\Omega)$.
The discrete weak form of~\eqref{eq:model} then reads: Find $u \in \SS^p_D(\TT_\coarse)$ such that
\begin{align}
\label{eq:discrete-system}
	&a(u_\coarse,v_\coarse)
	:=
	\int_\Omega \A \nabla u_\coarse \cdot \nabla v_\coarse + \vec{b} \cdot \nabla u_\coarse v_\coarse + c \, u_\coarse v_\coarse \d{x}
	+ \int_{\Gamma_R} \alpha \, u_\coarse v_\coarse \d{s} \\
\nonumber
	& \quad =
	\int_\Omega f v_\coarse + \vec{f} \cdot \nabla v_\coarse \d{x}
	+ \int_{\Gamma_N} \phi \, v_\coarse \d{s}
	+ \int_{\Gamma_R} \gamma \, v_\coarse \d{s}
	=:
	F(v_\coarse)
	\quad
	\text{for all } v_\coarse \in \SS^p_D(\TT_\coarse).
\end{align}

Our software is intended to solve equation~\eqref{eq:discrete-system} by adaptive finite element methods (AFEM), an abstract form of which is presented in the following~\cite{bangerth-rannacher,stevenson2007}.
\begin{algorithm}
\label{alg:afem}
	\textbf{Input}: Initial triangulation $\TT_0$ of $\Omega$\\
	\textbf{Loop}: For $\ell = 0, 1, \ldots$ do
	\begin{enumerate}[label={\rm (\roman*)}]
		\item Solve equation~\eqref{eq:discrete-system} to obtain $u_\ell$
		\item Estimate the error by computing refinement indicators $\eta_\ell(T)$ for all $T \in \TT_\ell$ 
		\item Mark elements $\MM_\ell \subseteq \TT_\ell$ based on $\eta_\ell$
		\item Refine marked elements to obtain $\TT_{\ell+1} := \refine(\TT_\ell, \MM_\ell)$
	\end{enumerate}
	\textbf{Output}: Sequence of solutions $u_\ell$
\end{algorithm}

A realization of the abstract adaptive Algorithm~\ref{alg:afem} is shown in the code snippet in Listing~\ref{lst:afem} below.
It computes the lowest-order FEM solution of the Poisson equation $-\Delta u = 1$ on the unit square $\Omega := (0,1)^2$ with homogeneous Dirichlet data $u=0$ on the boundary $\Gamma_D := \partial \Omega$.

\begin{lstlisting}[caption={Adaptive P1-FEM for Poisson problem $-\Delta u = 1$ on $\Omega = (0,1)^2$ subject to $u=0$ on $\partial \Omega$.},label={lst:afem}]
mesh = Mesh.loadFromGeometry('unitsquare');
fes = FeSpace(mesh, LowestOrderH1Fe);
u = FeFunction(fes);
blf = BilinearForm(fes);
blf.a = Constant(mesh, 1);
lf = LinearForm(fes);
lf.f = Constant(mesh, 1);

while mesh.nElements < 1e6
	A = assemble(blf);
	F = assemble(lf);
	freeDofs = getFreeDofs(fes);
	u.setFreeData(A(freeDofs,freeDofs) \ F(freeDofs));
	
	hT = sqrt(getAffineTransformation(mesh).area);
	qrEdge = QuadratureRule.ofOrder(1, '1D');
	qrTri = QuadratureRule.ofOrder(1, '2D');
	volumeRes = integrateElement(CompositeFunction(@(x) x.^2, lf.f), qrTri);
	edgeRes = integrateNormalJump(Gradient(u), qrEdge, @(j) j.^2, {}, ':');
	edgeRes(mesh.boundaries{:}) = 0;
	eta2 = hT.^2.*volumeRes + hT.*sum(edgeRes(mesh.element2edges), Dim.Vector);
	
	marked = markDoerflerSorting(eta2, 0.5);
	mesh.refineLocally(marked, 'NVB');
end
\end{lstlisting}

In Listing~\ref{lst:afem}, lines~1--7 are setup code that initializes all necessary data structures.
Lines~10--13 solve equation~\eqref{eq:discrete-system}, i.e., they realize Algorithm~\ref{alg:afem}{\rm (i)}.
For Algorithm~\ref{alg:afem}{\rm (ii)}, the error is then estimated in lines~15-21 by local contributions of the residual \textsl{a posteriori} error estimator~\cite{verfuerth}
\begin{equation*}
	\eta_\ell(T)^2
	:=
	|T| \, \norm{f}_{L^2(T)}^2 + |T|^{1/2} \norm{\jump{\nabla u \cdot \normalvec}}_{L^2(\partial T \cap \Omega)}^2
	\quad
	\text{for all } T \in \TT_\ell.
\end{equation*}
In line~23,  Algorithm~\ref{alg:afem}{\rm (iii)} is executed by the so-called \emph{Dörfler} marking criterion~\cite{doerfler1996}
\begin{equation*}
	\theta \sum_{T \in \TT_\ell} \eta_\ell(T)^2
	\leq
	\sum_{T \in \MM_\ell} \eta_\ell(T)^2
	\quad
	\text{with } \theta = 0.5.
\end{equation*}
Finally, line~24 corresponds to Algorithm~\ref{alg:afem}{\rm (iv)} and uses \emph{newest vertex bisection} (NVB)~\cite{stevenson2008} to refine (at least) the marked elements.

In the following, we assume that there is a fixed initial triangulation $\TT_0$ of $\Omega$.
All further meshes are supposed to be obtained by a finite number of successive refinement steps (with possibly multiple refinement strategies) from $\TT_0$.

\begin{remark}
\label{rem:blf-lf}
	We note that all other coefficients from~\eqref{eq:model} can be set just as easily as in line 6--7 of Listing~\ref{lst:afem}.
	The corresponding members of \lstinline|blf| and \lstinline|lf| are the following:
	\begin{itemize}
		\item \lstinline|blf.a|, \lstinline|blf.b|, \lstinline|blf.c|, \lstinline|lf.f|, and \lstinline|lf.fvec| for the data of~\eqref{eq:model-pde};
		
		\item \lstinline|blf.robin| and \lstinline|lf.robin| for the data of~\eqref{eq:model-robin};
		
		\item \lstinline|lf.neumann| for the data of~\eqref{eq:model-neumann}.
	\end{itemize}
	The types of functions that can be used are described in Section~\ref{sec:code-structure} below.
	There, also quadrature rules are presented, which can be assigned to the corresponding members \lstinline|qra|, \lstinline|qrb|, \lstinline|qrc|, and \lstinline|qrRobin| for \lstinline|blf| as well as \lstinline|qrf|, \lstinline|qrfvec|, \lstinline|qrRobin|, and \lstinline|qrNeumann| for \lstinline|lf|.
\end{remark}

\subsection{Necessity of OOP in \matlab FEM}

The use of OOP is not mandatory in \matlab, but it facilitates code that is powerful yet concise and flexible.
In particular, the code from Listing~\ref{lst:afem} relies heavily on OOP due to some peculiarities of the \matlab programming language.
Most notably, the referencing mechanics of \matlab differ greatly from languages like C/\cpp, Java, or Python, where referencing variables is either done by default, or explicitly.
We proceed by outlining the most important issues that we aim to address by OOP.

\subsubsection{Re-using data}\label{subsec:passing}

\matlab has a lazy copy policy for function arguments\additionalLink{\url{https://de.mathworks.com/help/matlab/matlab_prog/avoid-unnecessary-copies-of-data.html}}:
Arguments are generally passed by reference, but copied if they are modified within the function.
However, no local copy is made of variables that are assigned to themselves by returning data:

\begin{lstlisting}
function z = foo(x, y, z)
	y(1) = x;
	z = y + z;
end
z = foo(x, y, z)
\end{lstlisting}
In this example, \lstinline|x| is passed by reference, \lstinline|y| is copied because it is modified in line~2 (which is equivalent to passing the argument by value), and \lstinline|z| is modified but not copied since the output value is again assigned to \lstinline|z| in line~5.

FEM computations often re-use data throughout many sub-tasks; e.g., element areas are used in assembly of FEM systems, \textsl{a posteriori} error estimation, and even interpolation, via integration routines.
With the passing mechanics outlined above, there are two options to approach this issue:
First, data can be recomputed in each of the sub-tasks, yielding a clear public API, but causing superfluous operations.
Second, the data can be precomputed explicitly and held in memory.
In this case, the data management has to be done on the highest level of code by the user, or computations have to be grouped to respect data availability; both lead to a public API that is error-prone and inflexible.
It is therefore highly desirable to have proper call by reference mechanics, which, in \matlab, are only available through OOP.

\subsubsection{Encapsulating data efficiently}\label{subsec:value-handle}

The default for objects in \matlab are \emph{value classes}, which cannot change their state\additionalLink{\url{https://de.mathworks.com/help/matlab/matlab_oop/comparing-handle-and-value-classes.html}}.
The reason for this is the argument passing mechanism described above.
In fact, except for pathological cases, the method invocation \lstinline|obj.method(...)| and the function call \lstinline|method(obj, ...)| are equivalent\additionalLink{\url{https://de.mathworks.com/help/matlab/matlab_oop/method-invocation.html}}.
Hence, the copy-on-modify mechanics for function arguments also apply to instances of (value) classes.
However, classes that are derived from the abstract \lstinline|handle| class can overcome the limitations of value classes in the sense that they allow for modifications of state through methods:

\begin{minipage}{0.50\textwidth}
\begin{lstlisting}
classdef MyClass
	% ... 
	function obj = modify(obj, value)
		obj.field = value;
	end
end
obj = MyClass();
obj = modify(obj, 1);
\end{lstlisting}
\end{minipage}
\begin{minipage}{0.50\textwidth}
\begin{lstlisting}
classdef MyClass < handle
	% ...
	function modify(obj, value)
		obj.field = value
	end
end
obj = MyClass()
obj.modify(1)
\end{lstlisting}
\end{minipage}
Both code snippets result in \lstinline|obj.field| being equal to one.

Also, handle classes can inherently be referenced: assigning an instance of a handle class to a variable does not copy the underlying data, but only assigns a reference.
Finally, handle classes have native support for the \emph{observer pattern}, which is one of the central design elements of our code; see Section~\ref{subsec:observer} for details.

In order to communicate clearly, where methods can possibly alter the state of an object, we adhere to a coding best-practice called \emph{command query separation} throughout documentation and examples:
Commands, i.e., methods that alter the internal state of the calling object, are called with dot-notation \lstinline|obj.method(...)| and never return data;
queries, i.e., methods that do not alter the state of the calling object but may return data, are called with function call-notation \lstinline|method(obj, ...)|.

\subsubsection{Code safety and error checks}\label{subsec:statically-typed}

One of the disadvantages of dynamically typed languages like \matlab is the lack of automatic type checks and function overloading.
By using classes to represent (even trivial) data, this behavior can be emulated to a certain degree\additionalLink{\url{https://martinfowler.com/bliki/ValueObject.html}}.

Type checks can be automated by function argument validation, which was introduced recently in \matlab.
This is achieved by an optional \lstinline|arguments|-block after the function head that performs some checks on all input arguments of that function.
In particular, it can check class, dimension, and values of the input.
This greatly improves usability and error mitigation.

Since method invocation \lstinline|obj.method(...)| and function call \lstinline|method(obj, ...)| are virtually equivalent, function dispatch in \matlab goes by the first argument of a function.
This emulates function overloading at least in the first argument;
e.g., in Listing~\ref{lst:afem} (solving the Poisson problem), one can readily use \lstinline|plot(mesh)| and \lstinline|plot(u)| to plot the mesh and the FEM solution.
While this might be seen as syntactic sugar, it also greatly aids debugging.

\subsubsection{Vectorization}\label{subsec:vectorization}

One of the key features of \matlab are efficiently vectorized built-in linear algebra operations.
The usual local FEM formulation (i.e., on single elements and edges), however, does not allow for performance improvements through vectorization and parallelization, which are most pronounced if used with sufficiently large arrays to compensate for possible overhead.
It is therefore desirable to defer actual computation as long as possible to make optimal use of \matlab built-in routines.

Our code provides several well-defined interfaces very close to the natural (local) formulation of FEM which can be used to extend the functionality; see Section~\ref{subsec:module-fem}.
The referencing mechanics of handle classes then allow the internal generation of global data from this local code by pre-existing routines and passing it to efficient built-in routines.


\section{Code structure}\label{sec:code-structure}

\begin{figure}
	\centering
	\resizebox{0.9\textwidth}{!}{\begin{tikzpicture}

\newcommand{\intermediate}[2]{($(#1)!0.5!(#2)$)}
\pgfmathsetmacro{\pad}{0.3pt}
\pgfmathsetmacro{\ydist}{0pt}

\pgfdeclarelayer{background}
\pgfsetlayers{background,main}

\tikzstyle{module}[lightgray]=[draw,color=#1,fill=#1!50,very thick,text=#1!50!black,align=center,rounded corners,minimum width=11em,font=\bfseries]
\tikzstyle{class}=[align=left,minimum width=8em,text width=7.8em,font=\ttfamily\footnotesize\strut]
\tikzstyle{connection}[lightgray]=[color=#1,very thick,-{Diamond[length=0.7em,width=0.4em]}]
\tikzstyle{surround}[lightgray]=[draw,color=#1,fill=#1!20,rounded corners,very thick]

\node[module=color1] (M1) at (0,0) {Geometry};

\node[class,below=\pad of M1,xshift=0.5em] (C11) {Mesh};
\node[class,below=\ydist of C11] (C12) {AffineTransformation};
\node[class,below=\ydist of C12] (C13) {NVB};
\node[class,below=\ydist of C13] (C14) {AbstractBisection};

\node[module=color5,right=of M1.east] (M2) {Integration};

\node[class,below=\pad of M2,xshift=0.5em] (C21) {Evaluable};
\node[class,below=\ydist of C21] (C22) {QuadratureRule};
\node[class,below=\ydist of C22] (C23) {Barycentric};

\node[module=color2,right=of M2.east] (M3) {FEM};

\node[class,below=\pad of M3,xshift=0.5em] (C31) {FeSpace};
\node[class,below=\ydist of C31] (C32) {FiniteElement};
\node[class,below=\ydist of C32] (C33) {BilinearForm};
\node[class,below=\ydist of C33] (C34) {LinearForm};
\node[class,below=\ydist of C34] (C35) {Prolongation};

\draw[connection] (M1.south west) |- (C12.west);
\draw[connection] (M1.south west) |- (C13.west);
\draw[connection] (M1.south west) |- (C14.west);
\draw[connection=color1] (M1.south west) |- (C11.west);

\draw[connection] (M2.south west) |- (C22.west);
\draw[connection] (M2.south west) |- (C23.west);
\draw[connection=color5] (M2.south west) |- (C21.west);

\draw[connection] (M3.south west) |- (C32.west);
\draw[connection] (M3.south west) |- (C33.west);
\draw[connection] (M3.south west) |- (C34.west);
\draw[connection] (M3.south west) |- (C35.west);
\draw[connection=color2] (M3.south west) |- (C31.west);

\begin{pgfonlayer}{background}
	\coordinate (B3u) at ($(M1.north west)+(-3*\pad,3*\pad)$);
	\coordinate (B3l) at ($(M3.east |- C34.south)+(\pad,-3*\pad)$);
	\draw[surround=color2] (B3u) rectangle (B3l);
	
	\coordinate (B2u) at ($(M1.north west)+(-2*\pad,2*\pad)$);
	\coordinate (B2l) at ($(M2.east |- C34.south)+(\pad,-2*\pad)$);
	\draw[surround=color5] (B2u) rectangle (B2l);
	
	\coordinate (B1u) at ($(M1.north west)+(-\pad,\pad)$);
	\coordinate (B1l) at ($(M1.east |- C34.south)+(\pad,-\pad)$);
	\draw[surround=color1] (B1u) rectangle (B1l);
\end{pgfonlayer}

\end{tikzpicture}}
	\caption{Overall structure of the presented software package.
	Shown are all classes of the software package, subdivided into three modules.
	The most important class of each submodule is at the top of the respective list.}
	\label{fig:code-structure}
\end{figure}
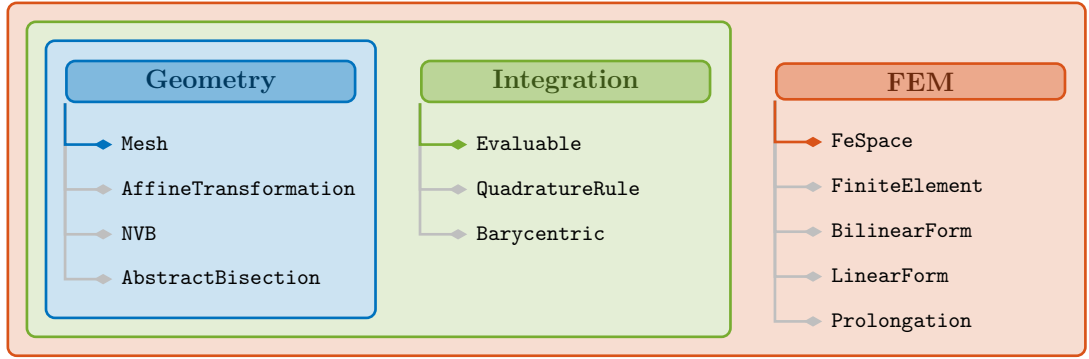

Most multi-purpose FEM packages have a huge code-base (often combining several languages) and, necessarily, a cleverly designed class hierarchy that may require significant effort to understand and get it running \cite{ngsolve,fenics,dealII,dune,bhl+2021}.
Since one of our aims is that our code is easy to modify and extend, we strive for a relatively small code-base that nevertheless covers the most widespread demands of academic FEM software.
Thus, the core of our software package is made up of only twelve classes and interfaces that can be roughly divided into three modules: Geometry, integration, and FEM.
This partition is shown in Figure~\ref{fig:code-structure}.

What follows is a description of the separate core classes as well as their inter relationship to one another.
Following the partition of the code, our presentation is divided into three parts.

\subsection{Module geometry}\label{subsec:module-geometry}

The geometry module can be used entirely on its own.
It handles mesh generation as well as local mesh refinement.

\subsubsection{Mesh representation}
As the underlying mesh $\TT_\coarse$ is the cornerstone in any adaptive FEM algorithm, the \lstinline|Mesh| class is the central building block of \moofem.
It consists of all data that is important for the digital representation of a 2D triangulation:
coordinates, edges, elements, connectivity of edges and elements, edge orientation, and boundary information.
The precise data structures for storing this information are described in Section~\ref{sec:data-structures} below.

Here, we focus on the role of the class in the code compound.
Most other classes need a \lstinline|Mesh| to function properly and, hence, store a reference to a suitable instance.
Validity of data is strongly coupled to the underlying mesh: as soon as the mesh changes, derived data (e.g., geometric information as well as the data used in FEM computations) may be invalid.
It is therefore of vital importance that changes in the mesh are made public to every object that stores a reference to it, contrary to the usual flow of information in object oriented code.
Also, objects that do not explicitly depend on a mesh may need to act if the mesh changes.
A well-known remedy for this issue is presented in the next section.

\subsubsection{Observer pattern}\label{subsec:observer}
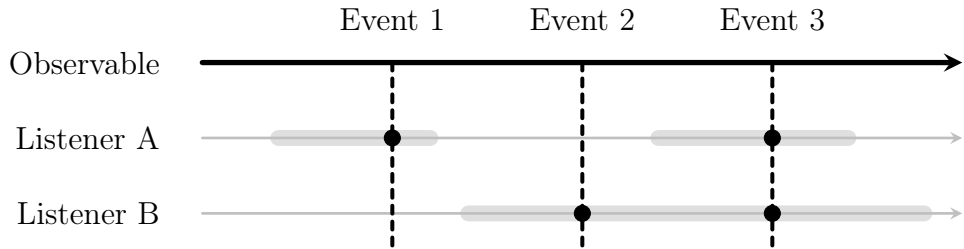
\begin{figure}
	\centering
	\begin{tikzpicture}
\pgfmathsetmacro{\xend}{10}
\pgfmathsetmacro{\ydist}{-1}

\tikzstyle{lifetime}=[-stealth,line width=1pt,line cap=round]
\tikzstyle{event}=[dashed,line width=1.5pt,line cap=round]
\tikzstyle{subscription}=[line width=6pt,line cap=round,lightgray,draw opacity=0.5]
\tikzstyle{reaction}=[draw,fill,circle,minimum size=6pt,inner sep=0pt]

\coordinate (OB) at (0,0);
\coordinate (L1) at (0,\ydist);
\coordinate (L2) at (0,2*\ydist);

\coordinate (E1) at (0.25*\xend,0);
\coordinate (E2) at (0.50*\xend,0);
\coordinate (E3) at (0.75*\xend,0);

\node[left=1em of OB] {Observable};
\node[left=1em of L1] {Listener A};
\node[left=1em of L2] {Listener B};

\node[above=0.7em of E1] {Event 1};
\node[above=0.7em of E2] {Event 2};
\node[above=0.7em of E3] {Event 3};

\draw[lifetime,line width=2pt] (OB) -- +(\xend,0);
\draw[lifetime,lightgray] (L1) -- +(\xend,0);
\draw[lifetime,lightgray] (L2) -- +(\xend,0);

\draw[subscription] ($(L1)+(0.1*\xend,0)$) to +(0.2*\xend,0);
\draw[subscription] ($(L1)+(0.6*\xend,0)$) to +(0.25*\xend,0);
\draw[subscription] ($(L2)+(0.35*\xend,0)$) to +(0.6*\xend,0);

\draw[event] (E1) -- +(0,2.5*\ydist);
\draw[event] (E2) -- +(0,2.5*\ydist);
\draw[event] (E3) -- +(0,2.5*\ydist);

\node[reaction] at (0.25*\xend,\ydist) {};
\node[reaction] at (0.75*\xend,\ydist) {};
\node[reaction] at (0.50*\xend,2*\ydist) {};
\node[reaction] at (0.75*\xend,2*\ydist) {};

\end{tikzpicture}
	\caption{Schematic functionality of the observer pattern.
	The timelines represent the lifetimes of the observable object and the listeners.
	Two listeners are temporarily registered to receive events (highlighted in gray).
	If the listeners are registered during the broadcast of an event (dashed lines), some internal reaction is triggered (bold dots).
	For each additional observable or event type of the same observable, a separate graph is needed.}
	\label{fig:observer}
\end{figure}

The observer pattern~\cite{design-patterns} is used to broadcast \emph{events} from a central object (the \emph{observable}) to other objects (the observers, termed \emph{listeners} in \matlab), which can react to the event in a predefined way and which can, at runtime, register and de-register to receive such events\additionalLink{\url{https://de.mathworks.com/help/matlab/matlab_oop/learning-to-use-events-and-listeners.html}}; see Figure~\ref{fig:observer}.
The events signal, e.g., a change of state of the observable.
All classes that derive from \lstinline|handle| in \matlab automatically implement the interfaces necessary to be the source of events.
In particular, our \lstinline|Mesh| class broadcasts events that signal a change in the mesh (e.g., by the call $\refine(\TT_\ell, \MM_\ell)$ in Algorithm~\ref{alg:afem}) as well as completed computation of bisection data to signal imminent refinement; this is covered in the next two sections.

\subsubsection{Refinement}

\begin{figure}
	\centering
	\begin{tikzpicture}

\tikzstyle{triangle}=[thick]

\newcommand{\referenceEdge}[3]{%
	\coordinate (M) at ($0.333*#1+0.333*#2+0.333*#3$);%
	\coordinate (Am) at ($#1!0.3!(M)$);%
	\coordinate (Bm) at ($#2!0.3!(M)$);%
	\draw[triangle,semithick,red] (Am) -- (Bm);}

\newcommand{\drawParentTriangle}{%
	\coordinate (P1) at (0,0);
	\coordinate (P2) at (3,0);
	\coordinate (P3) at (2.1,1.5);
	\draw[triangle] (P1) -- (P2) -- (P3) -- cycle;
}

\pgfmathsetmacro{\xdist}{3.5}
\pgfmathsetmacro{\ydist}{-2.0}

\begin{scope}
	\drawParentTriangle
	\coordinate (N1) at ($0.5*(P1)+0.5*(P2)$);
	\draw[triangle] (P3) -- (N1);
	\referenceEdge{(P3)}{(P1)}{(N1)}
	\referenceEdge{(P2)}{(P3)}{(N1)}
\end{scope}

\begin{scope}[shift={(1*\xdist,0)}]
	\drawParentTriangle
	\coordinate (N1) at ($0.5*(P1)+0.5*(P2)$);
	\coordinate (N2) at ($0.5*(P2)+0.5*(P3)$);
	\draw[triangle] (P3) -- (N1);
	\draw[triangle] (N1) -- (N2);
	\referenceEdge{(P3)}{(P1)}{(N1)}
	\referenceEdge{(N1)}{(P2)}{(N2)}
	\referenceEdge{(P3)}{(N1)}{(N2)}
\end{scope}

\begin{scope}[shift={(2*\xdist,0)}]
	\drawParentTriangle
	\coordinate (N1) at ($0.5*(P1)+0.5*(P2)$);
	\coordinate (N3) at ($0.5*(P1)+0.5*(P3)$);
	\draw[triangle] (P3) -- (N1);
	\draw[triangle] (N1) -- (N3);
	\referenceEdge{(P2)}{(P3)}{(N1)}
	\referenceEdge{(P1)}{(N1)}{(N3)}
	\referenceEdge{(N1)}{(P3)}{(N3)}
\end{scope}

\begin{scope}[shift={(0*\xdist,1*\ydist)}]
	\drawParentTriangle
	\coordinate (N1) at ($0.5*(P1)+0.5*(P2)$);
	\coordinate (N2) at ($0.5*(P2)+0.5*(P3)$);
	\coordinate (N3) at ($0.5*(P1)+0.5*(P3)$);
	\draw[triangle] (P3) -- (N1);
	\draw[triangle] (N3) -- (N1);
	\draw[triangle] (N1) -- (N2);
	\referenceEdge{(P1)}{(N1)}{(N3)}
	\referenceEdge{(N1)}{(P2)}{(N2)}
	\referenceEdge{(P3)}{(N1)}{(N3)}
	\referenceEdge{(N1)}{(P3)}{(N2)}
\end{scope}

\begin{scope}[shift={(1*\xdist,1*\ydist)}]
	\drawParentTriangle
	\coordinate (N1) at ($0.5*(P1)+0.5*(P2)$);
	\coordinate (N2) at ($0.5*(P2)+0.5*(P3)$);
	\coordinate (N3) at ($0.5*(P1)+0.5*(P3)$);
	\coordinate (N4) at ($0.333*(P1)+0.333*(P2)+0.333*(P3)$);
	\draw[triangle] (P3) -- (N1);
	\draw[triangle] (N3) -- (N1);
	\draw[triangle] (N1) -- (N2);
	\draw[triangle] (N2) -- (N3);
	\referenceEdge{(P1)}{(N1)}{(N3)}
	\referenceEdge{(N1)}{(P2)}{(N2)}
	\referenceEdge{(N1)}{(N2)}{(N4)}
	\referenceEdge{(N2)}{(P3)}{(N4)}
	\referenceEdge{(P3)}{(N3)}{(N4)}
	\referenceEdge{(N3)}{(N1)}{(N4)}
\end{scope}

\begin{scope}[shift={(2*\xdist,1*\ydist)}]
	\drawParentTriangle
	\coordinate (N1) at ($0.5*(P1)+0.5*(P2)$);
	\coordinate (N2) at ($0.5*(P2)+0.5*(P3)$);
	\coordinate (N3) at ($0.5*(P1)+0.5*(P3)$);
	\draw[triangle] (N3) -- (N1);
	\draw[triangle] (N1) -- (N2);
	\draw[triangle] (N2) -- (N3);
	\referenceEdge{(P1)}{(N1)}{(N3)}
	\referenceEdge{(N1)}{(P2)}{(N2)}
	\referenceEdge{(N2)}{(N3)}{(N1)}
	\referenceEdge{(N3)}{(N2)}{(P3)}
\end{scope}

\end{tikzpicture}
	\caption{Implemented bisection methods (top left to bottom right): \lstinline|Bisec1|, \lstinline|Bisec12|, \lstinline|Bisec13|, \lstinline|Bisec123|, \lstinline|Bisec5|, and \lstinline|BisecRed|.
	The refinement edge of the parent triangle is the bottom line, those of the children are highlighted by parallel lines.}
	\label{fig:bisections}
\end{figure}
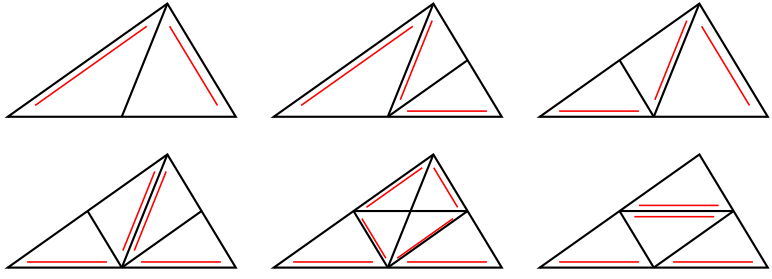

Mesh refinement is implementationally divided into bisection of single elements $T \in \TT_\coarse$ and coordination of bisections on the whole mesh $\TT_\coarse$ to obtain $\TT_\fine = \refine(\TT_\coarse, \MM_\coarse)$.
In the class \lstinline|AbstractBisection|, possible bisections of a single triangle $T \in \TT_\coarse$ are encoded.
Subclasses of this class manage the generation of all \lstinline|Mesh| data structures (see Section~\ref{sec:data-structures} below) in the passage from $T \in \TT_\coarse$ to its children $\set{T' \in \TT_\fine}{T' \subseteq T}$  in $\TT_\fine$.
The subclasses that are currently implemented are shown in Figure~\ref{fig:bisections} and build on the implementation of~\cite{p1afem}.

Local bisections are combined in mesh refinement routines derived from newest vertex bisection (NVB)~\cite{stevenson2008,kpp2013}, in the course of which all elements $T \in \TT_\coarse$ are assigned a subclass of \lstinline|AbstractBisection|.
We follow the iterative algorithm given in~\cite{kpp2013} for NVB, which terminates regardless of the mesh $\TT_\coarse$ under consideration.
Denoting the set of edges in $\TT_\coarse$ by $\EE_\coarse$ and given a subset $\MM_\coarse \subseteq \TT_\coarse$ of marked elements, the abstract scheme is the following:
\begin{enumerate}[label={\rm (\roman*)}]
	\item Determine all edges in $\EE_\coarse$ that have to be bisected in order to bisect all elements in $\MM_\coarse$ by a given bisection rule.
	\item Compute the \emph{mesh closure}, i.e., determine all edges that \emph{additionally} have to be refined to eliminate hanging vertices.
	\item For every element in $\TT_\coarse$ determine which bisection method to employ, based on the marked edges.
	\item Execute bisection of all elements according to their assigned bisection methods.
\end{enumerate}
The first three steps are executed by subclasses of \lstinline|NVB|, the fourth step is carried out by the mesh itself.
The rationale behind this splitting is that, before the fourth step, all necessary information to carry out mesh refinement is already known; thus, this provides a natural hook for other classes to harness this information, e.g., prolongation operators or multi-grid solvers.

Several realizations of the abstract scheme presented above are implemented in our software package: \lstinline|NVB1|, \lstinline|NVB3| (= \lstinline|NVB|), \lstinline|NVB5|, and \lstinline|RGB|, which are outlined in~\cite{p1afem}, as well as \lstinline|NVBEdge|, which is an edge-driven refinement strategy used in~\cite{dks2016,ip2021}.

\subsubsection{Additional geometric information}

For triangular meshes, all elements can be obtained by affine transformations of the so-called reference triangle $T_\mathrm{ref} := \conv\{ (0,0), \allowbreak (1,0), (0,1) \}$, i.e., for all $T \in \TT_\coarse$, there exists an affine diffeomorphism $F_T \colon T_\mathrm{ref} \to T$.
Many computations in FEM need additional geometric information based on this diffeomorphism.
In particular, the transposed inverse $\D F_T^{-\intercal}$ and the determinant $\det\D F_T = 2 |T|$ of its derivative are of utmost importance in integration and assembly routines.
Together with the length $\d{s} = |E|$ of each edge $E \in \EE_\coarse$ and its unit-normal vector $\normalvec_E$, these data are stored in instances of \lstinline|AffineTransformation|, which take a \lstinline|Mesh| to construct.

For performance reasons, instances of \lstinline|AffineTransformation| are requested from the mesh and are cached, i.e., computed at the first request, then stored as a reference within the mesh object and returned if further requests occur.

\subsection{Module integration}\label{subsec:module-integration}

Integration is a crucial part of FEM assembly and postprocessing (e.g., \textsl{a posteriori} error estimation).
The module defines two classes that encapsulate data for numerical integration (also termed \emph{quadrature}) routines and one to encapsulate function evaluation with a unified interface.
This module can be used only in conjunction with the geometry module.

\subsubsection{Barycentric coordinates}\label{subsec:barycentric}

We denote all evaluation points and quadrature nodes in barycentric coordinates.
On a triangle $T = \conv\{ z^{(1)}, z^{(2)}, z^{(3)} \}$ and a point $x \in T$, we denote by $\lambda \in [0,1]^3$ the barycentric coordinates of $x$, determined by the equations
\begin{equation*}
\sum_{i=1}^{3} \lambda_i = 1
\quad \text{and} \quad
\sum_{i=1}^{3} \lambda_i z^{(i)} = x.
\end{equation*}
If the triangle $T$ is non-degenerate, $\lambda$ is unique.
In \moofem, 2D barycentric coordinates are implemented in the class \lstinline|Barycentric2D|, which is derived from the abstract class \lstinline|Barycentric|.
For convenience, this class stores a collection of barycentric coordinates.

The concept of barycentric coordinates can be generalized to $d$-simplices with $d \geq 1$ such that $\lambda \in [0,1]^{d+1}$.
In particular, in the case $d=1$ we have $\lambda \in [0,1]^2$, which is used for evaluation points and quadrature rules on edges and implemented in the class \lstinline|Barycentric1D|.
For any point $x \in T_\mathrm{ref} = \conv\{ (0,0), (1,0), (0,1) \}$ in the reference triangle, the barycentric coordinates are $\lambda = (1-x_1-x_2,x_1,x_2)$;
for any point $x \in E_\mathrm{ref} = [0,1]$ on the reference edge, they are $\lambda = (x, 1-x)$.

Denoting all function evaluation points in barycentric coordinates allows for triangle independent representation.
This is reflected by \lstinline|Evaluable.eval| below, which is the core of our vectorization efforts and paramount for the efficiency of \moofem.

\subsubsection{Quadrature data}

To approximate the integral of a (possibly vector-valued) function $f \colon T \to \R^n$ for some $n \in \N$ on a triangle $T \in \TT_\coarse$, we employ numerical quadrature:
\begin{equation}
\label{eq:quadrature-element}
	\int_T f(x) \d{x}
	\approx
	|T| \sum_{k=1}^{N} \omega_k f\big( x(\lambda^{(k)}, T) \big),
\end{equation}
where $(\lambda^{(k)})_{k=1}^N$ is a collection of barycentric coordinates, $x(\lambda^{(k)}, T)$ is the Cartesian coordinate corresponding to the $k$-th barycentric coordinate on $T$, and $(\omega_k)_{k=1}^N$ are weights with $\sum_{k=1}^N \omega_k = 1$.

Barycentric coordinates and weights make up a \lstinline|QuadratureRule| object.
Quadrature rules can either be constructed explicitly by giving barycentric coordinates and weights, or by the static method
\begin{lstlisting}
qr = QuadratureRule.ofOrder(order, [dim]);
\end{lstlisting}
The optional string argument \lstinline|dim| is used to distinguish between \lstinline|'1D'| and \lstinline|'2D'| quadrature rules, the default being \lstinline|'2D'|.
For 1D, suitable Gauss-rules are implemented.
For 2D, symmetric quadrature rules up to order 5 are implemented~\cite{zcl2009}.
Higher order quadrature rules on triangles use tensorized Gauss-rules on $[0,1]^2$ and the Duffy transform $\Phi \colon [0,1]^2 \to T_\mathrm{ref}, \Phi(s,t) = \big(s, t(1-s)\big)$;~\cite{duffy1982}.

\subsubsection{Function evaluation}\label{subsec:evaluable}

The core of the integration module is a wrapper for functions $f \colon \Omega \to \R^n$ for some $n \in \N$, the abstract \lstinline|Evaluable| class:
\begin{lstlisting}
classdef Evaluable < handle
	properties (Abstract, GetAccess='public', SetAccess='protected')
		mesh
	end
	methods (Abstract, Access='public')
		eval(obj, bary, idx)
	end
end
\end{lstlisting}
The abstract method \lstinline|eval| is intended to evaluate the function at all points $x(\lambda^{(k)}, T)$ for all barycentric coordinates $\lambda^{(k)}$ in \lstinline|bary| and all elements given by the index set \lstinline|idx|.

By programming all routines (e.g., integration, plotting, finite element assembly) only to this interface, one can wrap virtually anything in a suitable subclass of \lstinline|Evaluable| and readily use the predefined routines.
The most important classes that implement this interface are:
\begin{itemize}
	\item \lstinline|Constant|: Efficiently wraps constant functions in the \lstinline|Evaluable| interface.
	
	\item \lstinline|MeshFunction|: General functions $f \colon \Omega \to \R^n$ for some $n \in \N$.
	
	\item \lstinline|FeFunction|: FEM functions, e.g., $u \in \SS^p(\TT_\coarse)$.
	
	\item \lstinline|Gradient|, \lstinline|Hessian|: Element-wise gradient $\nabla u$ and Hessian $\nabla^2 u$ for FEM functions.
	
	\item \lstinline|CompositeFunction|: Arbitrarily combine any of the above; see the following explanation.
\end{itemize}
In particular, all of the above can be used as coefficients in~\eqref{eq:model}.

Especially powerful is the subclass \lstinline|CompositeFunction|, which uses the composite pattern~\cite{design-patterns}:
\begin{lstlisting}
f = CompositeFunction(funcHandle, funcArgument1, ..., funcArgumentN)
\end{lstlisting}
This class takes a function handle and one \lstinline|Evaluable| for every argument of the function handle.
E.g., $x u^2$ can be implemented by
\begin{lstlisting}
f = CompositeFunction(@(x,u) x.*u.^2, ...
			MeshFunction(mesh, @(x) x(1,:,:)), FeFunction(fes));
\end{lstlisting}
If \lstinline|f| is evaluated, first all arguments are evaluated, then the resulting data is processed by the function handle.
By polymorphism, the \lstinline|Evaluable| arguments can be of any subclass.
This allows to define complex functions that still can be evaluated efficiently, since evaluation of the function handle is deferred until the data for all requested elements and quadrature nodes is available.
In this way, the vectorization capabilities of \matlab are used to full extent.

\subsubsection{Quadrature routines}

There are several routines for quadrature implemented in our \moofem package:
\begin{itemize}
	\item \lstinline|integrateElement|: Integration on elements; cf.~\eqref{eq:quadrature-element}.
	
	\item \lstinline|integrateEdge|: An analogue over edges.
	This can only be used for subclasses of \lstinline|Evaluable| that implement the method \lstinline|evalEdge|.
	Edge evaluation is not well-defined for some functions that are not continuous across edges (e.g., element-wise polynomials).
	
	\item \lstinline|integrateJump|, \lstinline|integrateNormalJump|: Integrate $\jump{\cdot}$ and $\jump{(\cdot) \cdot \normalvec}$ over edges, respectively.
\end{itemize}
All quadrature routines take an \lstinline|Evaluable|, a \lstinline|QuadratureRule|, and an optional set of indices that corresponds to a subset of elements or edges on which the quadrature should be evaluated.
The routines handling (normal) jumps take additional arguments: a function handle, a cell array of \lstinline|Evaluable|, and edge indices.
\begin{lstlisting}
int = integrateJump(f, qr, funcHandle, funcArg, idx)
\end{lstlisting}
This is used for post-processing the jump with the first argument of the function handle being reserved for the jump; i.e., the routine works roughly as follows:
First, compute the jump by $\texttt{jump} = \jump{\texttt{f}}$ (or $\texttt{jump} = \jump{\texttt{f}\cdot \normalvec}$).
Then, all additional \lstinline|Evaluables| are evaluated on the edges indicated by \lstinline|idx|, where also the value of the jump is updated by the function handle:
\begin{lstlisting}
val = {evalEdge(funcArg{1}, idx), ..., evalEdge(funcArg{N}, idx)};
jump(idx) = functionHandle(jump(idx), val{1}, ..., val{N});
\end{lstlisting}
Multiple such triplets \lstinline|funcHandle|, \lstinline|funcArg|, \lstinline|idx| for post-processing are allowed and sequentially applied as in the above listing, one after another.

\subsection{Module FEM}\label{subsec:module-fem}

This last module uses the classes presented in the last sections to conveniently represent FEM functions and efficiently assemble FEM data.

\subsubsection{Finite element spaces}

Finite elements are usually defined on the reference triangle $T_\mathrm{ref}$, everything else follows from using the affine transformation $F_T$ for every $T \in \TT_\coarse$.
This viewpoint is represented accordingly in our code.
The abstract class \lstinline|FiniteElement| asks to implement evaluation of all basis functions (as well as their gradient and their Hessian) on $T_\mathrm{ref}$.
Furthermore, evaluation on a reference edge (if applicable) and the connectivity of the degrees of freedom (DOFs), i.e., how the finite element couples across edges and vertices, have to be specified.

The class \lstinline|FeSpace| combines the local information of finite elements with the global geometry of the mesh.
It takes a \lstinline|Mesh| and a \lstinline|FiniteElement| to assemble lists of global DOFs per element and edge, as well as free DOFs, i.e., DOFs that do not lie on $\Gamma_D$.

So far, Lagrange finite elements of arbitrary order are implemented; both $H^1(\Omega)$-conforming (i.e., $\SS^p(\TT_\coarse) =  \PP^p(\TT_\coarse) \cap H^1(\Omega)$) and $L^2(\Omega)$-conforming (i.e., $\PP^p(\TT_\coarse)$).
The underlying implementation uses Bernstein--B\'ezier polynomials~\cite{aad2011}.
For lowest-order finite elements (both continuous and discontinuous), additional optimized implementations are available.

\subsubsection{FEM system assembly}

Let $(\varphi_k)_{k=1}^N$ be a basis of $\SS^p(\TT_\coarse)$.
Responsible for the assembly of the data
\begin{subequations}
\label{eq:fem-data}
\begin{align}
	A_{ij}
	&:=
	\int_\Omega \A \nabla \varphi_j \cdot \nabla \varphi_i + \vec{b} \cdot \nabla \varphi_j \varphi_i + c \, \varphi_j \varphi_i \d{x}
	+ \int_{\Gamma_R} \alpha \, \varphi_j \varphi_i \d{s},\\
	F_i
	&:=
	\int_\Omega f \varphi_i + \vec{f} \cdot \nabla \varphi_i \d{x}
	+ \int_{\Gamma_N} \phi \, \varphi_i \d{s}
	+ \int_{\Gamma_R} \gamma \, \varphi_i \d{s}
\end{align}
\end{subequations}
are the classes \lstinline|BilinearForm| and \lstinline|LinearForm|, respectively.
Both classes have fields for their respective coefficients and quadrature rules for each of the terms in~\eqref{eq:fem-data}; see Remark~\ref{rem:blf-lf}.
The data in~\eqref{eq:fem-data} are then obtained by calling the \lstinline|assemble| methods of both classes, using a generalization of the efficiently vectorized method outlined in~\cite{p1afem}.

Note that \lstinline|Evaluable| is a handle class.
Thus, no data must be copied to set (bi-)linear form coefficients.
Furthermore, in situations where the coefficients change frequently, e.g., in the presence of iterative solvers for nonlinear PDEs, the coefficients of the (bi-)linear form can change between two consecutive calls of \lstinline|assemble|.
This results in much cleaner code since one does not need to re-set the coefficients.

\subsubsection{Prolongation}

It is often necessary to prolongate FEM functions $u_\coarse \in \SS^p(\TT_\coarse)$ on a mesh $\TT_\coarse$ to the richer space $\SS^p(\TT_\fine)$ with a refined mesh $\TT_\fine$.
This is handled by subclasses of the abstract class \lstinline|Prolongation|.
The implemented prolongation methods are \lstinline|LoFeProlongation| and \lstinline|FeProlongation| for lowest-order and general (in particular, higher-order) FEM functions, respectively.
Note that the latter is not tailored to a specific finite element and, hence, its computational effort is slightly higher than that of the former.
The syntax of prolongation of a function \lstinline|u = FeFunction(fes)| on a finite element space \lstinline|fes| from a coarse to a refined mesh is as follows:
\begin{lstlisting}
P = FeProlongation(fes);
mesh.refineLocally(marked);
u.setData(P.prolongate(u));
\end{lstlisting}
The call to \lstinline|P.prolongate| performs a matrix-vector multiplication with the (sparse) prolongation matrix \lstinline|P.matrix|, which is set automatically by element-wise Lagrange interpolation whenever the mesh is refined, due to the events sent by the mesh; see Section~\ref{subsec:observer} and the examples in Section~\ref{sec:examples}.


\section{Data structures}\label{sec:data-structures}

\subsection{Mesh}
The \lstinline|Mesh| class stores information about coordinates, edges, and elements.
In the following, let $\nC, \nS, \nE \in \N$ denote the number of vertices, edges, and elements, respectively.
The class has the following fields:
\begin{itemize}
	\item \lstinline|coordinates| ($2 \times \nC$):
	Coordinates of mesh vertices.
	The entries \lstinline|coordinates(1,i)| and \lstinline|coordinates(2,i)| store the $x$- and $y$-coordinates of the $i$-th vertex, respectively.
	The order of the coordinates is provided by the user.
	
\item \lstinline|edges| ($2 \times \nS$):
	Indices of vertices of all edges in the mesh.
	The $i$-th edge of the mesh starts at vertex \lstinline|edges(1,i)| and ends at vertex \lstinline|edges(2,i)|.
	The order is determined automatically from information provided in \lstinline|elements|.
	Boundary edges are oriented such that the domain lies on its left; inner edges cannot be assigned a meaningful orientation and, therefore, it is chosen randomly.
	
\item \lstinline|elements| ($3 \times \nE$):
	Indices of vertices of which elements are comprised.
	The $i$-th element is spanned by the vertices with indices \lstinline|elements(:,i)|, where the order is counter-clockwise.
	The order of elements is provided by the user.
	
\item \lstinline|element2edges| ($3 \times \nE$):
	Indices of edges which are contained in elements.
	The $i$-th element contains edges with indices \lstinline|element2edges(:,i)|.
	The $j$-th edge \lstinline|element2edges(j,i)| of the $i$-th element is the one between the vertices with indices \lstinline|elements(j,i)| and \lstinline|elements(mod(j+1,3)+1, i)| (but not necessarily in that order).
	This information is determined automatically.
	
\item \lstinline|boundaries| (cell array):
	Indices of all edges that form a specific part of the boundary.
	The cell \lstinline|boundaries{i}| is a vector of edge indices that form the $i$-th boundary (if present).
	The boundary parts are provided by the user (see below), but the association with edge indices is done automatically.
\end{itemize}

The orientation of the normal vector from \lstinline|AffineTransformation| follows the orientation of the edges: it points to the right of the edge.
In particular, this means that the normal vectors on boundary edges point out of the domain.

\subsection{Mesh construction}
An instance of the \lstinline|Mesh| class can be constructed by
\begin{lstlisting}
mesh = Mesh(coordinates, elements, bndEdges);
\end{lstlisting}
Here, \lstinline|coordinates| and \lstinline|elements| have to be given as they are stored in the class (see above).
The cell array \lstinline|bndEdges| describes the boundary parts, where the $i$-th edge of the $k$-th boundary part lies between the vertices with indices \lstinline|bndEdges{k}(1,i)| and \lstinline|bndEdges{k}(2,i)|.
This is necessary because the user does not know the internal edge numbering before construction.
Special attention has to be payed to the correct orientation of the elements (counter-clockwise) and the edges on the boundary (domain on their left), because this is not checked by the constructor.

The arrays required by the constructor can be assembled and passed from a \matlab script, or loaded from comma separated value files by the static method
\begin{lstlisting}
mesh = Mesh.loadFromGeometry('<name>');
\end{lstlisting}
These files must be placed in a subdirectory \lstinline|lib/mesh/@Mesh/geometries/<name>| and be named \lstinline|coordinates.dat|, \lstinline|elements.dat|, and \lstinline|boundary<n>.dat|, where boundary parts are enumerated by $n \in \N$.

\begin{figure}
	\centering
	\resizebox{0.25\textwidth}{!}{
\begin{tikzpicture}[baseline=(current bounding box.center)]

\pgfmathsetmacro{\len}{5}

\tikzstyle{edge}=[line width=2pt,black]
\tikzstyle{boundary}=[line width=6pt,line cap=round,draw opacity=0.5,rounded corners=2pt]
\tikzstyle{coord}=[draw,fill,circle,minimum size=8pt,inner sep=0pt]
\tikzstyle{mylabel}=[outer sep=3pt]

\coordinate (P1) at (0,0);
\coordinate (P2) at (\len,0);
\coordinate (P3) at (\len,\len);
\coordinate (P4) at (0,\len);
\coordinate (P5) at (0.5*\len,0.5*\len);

\coordinate (M1) at ($0.333*(P1)+0.333*(P2)+0.333*(P5)$);
\coordinate (M2) at ($0.333*(P2)+0.333*(P3)+0.333*(P5)$);
\coordinate (M3) at ($0.333*(P3)+0.333*(P4)+0.333*(P5)$);
\coordinate (M4) at ($0.333*(P4)+0.333*(P1)+0.333*(P5)$);

\coordinate (L1) at ($0.5*(P1)+0.5*(P2)$);
\coordinate (L2) at ($0.5*(P1)+0.5*(P4)$);
\coordinate (L3) at ($0.5*(P1)+0.5*(P5)$);
\coordinate (L4) at ($0.5*(P2)+0.5*(P3)$);
\coordinate (L5) at ($0.5*(P2)+0.5*(P5)$);
\coordinate (L6) at ($0.5*(P3)+0.5*(P4)$);
\coordinate (L7) at ($0.5*(P3)+0.5*(P5)$);
\coordinate (L8) at ($0.5*(P4)+0.5*(P5)$);

\draw[edge] (P1) -- (P2) -- (P3) -- (P4) -- cycle;
\draw[edge] (P1) -- (P3);
\draw[edge] (P2) -- (P4);

\draw[boundary,color=color2] (P4) -- (P1) -- (P2);
\draw[boundary,color=color5] (P2) -- (P3) -- (P4);

\node[coord] at (P1) {};
\node[coord] at (P2) {};
\node[coord] at (P3) {};
\node[coord] at (P4) {};
\node[coord] at (P5) {};

\node[mylabel,anchor=north east] at (P1) {1};
\node[mylabel,anchor=north west] at (P2) {2};
\node[mylabel,anchor=south west] at (P3) {3};
\node[mylabel,anchor=south east] at (P4) {4};
\node[mylabel,anchor=south] at (P5) {5};

\node[mylabel] at (M1) {1};
\node[mylabel] at (M2) {2};
\node[mylabel] at (M3) {3};
\node[mylabel] at (M4) {4};

\node[mylabel,anchor=north] at (L1) {1};
\node[mylabel,anchor=east] at (L2) {2};
\node[mylabel,anchor=south] at (L3) {3};
\node[mylabel,anchor=west] at (L4) {4};
\node[mylabel,anchor=south] at (L5) {5};
\node[mylabel,anchor=south] at (L6) {6};
\node[mylabel,anchor=south] at (L7) {7};
\node[mylabel,anchor=south] at (L8) {8};

\end{tikzpicture}}
\quad%
\newcolumntype{A}[1]{ >{\raggedleft\arraybackslash} m{#1} }
\newcolumntype{B}{ >{\centering\arraybackslash} m{0.5cm} }
\resizebox{0.42\textwidth}{!}{
\begin{tabular}{A{2.8cm}|BBBBBBBB}
	$n$                                     & 1   & 2   & 3   & 4   & 5   & 6   & 7   & 8  \\
	\hline
	\multirow{2}{*}{\texttt{coordinates}}   & 0.0 & 1.0 & 1.0 & 0.0 & 0.5 &     &     &    \\
                                            & 0.0 & 0.0 & 1.0 & 1.0 & 0.5 &     &     &    \\
    \hline
	\multirow{2}{*}{\texttt{edges}}         & 1   & 4   & 1   & 2   & 2   & 3   & 3   & 4  \\
	                                        & 2   & 1   & 5   & 3   & 5   & 4   & 5   & 5  \\
	\hline
	\multirow{3}{*}{\texttt{elements}}      & 1   & 2   & 3   & 4   &     &     &     &    \\
	                                        & 2   & 3   & 4   & 1   &     &     &     &    \\
	                                        & 5   & 5   & 5   & 5   &     &     &     &    \\
	\hline
	\multirow{3}{*}{\texttt{element2edges}} & 1   & 4   & 6   & 2   &     &     &     &    \\
	                                        & 5   & 7   & 8   & 3   &     &     &     &    \\
	                                        & 3   & 5   & 7   & 8   &     &     &     &    \\
\end{tabular}}
\quad%
\resizebox{0.23\textwidth}{!}{
\begin{tabular}{A{2.8cm}|BB}
	$n$ & 1 & 2 \\
	\hline
	\multirow{2}{*}{\texttt{bndEdges\{1\}}} & 1 & 4 \\
	                                        & 2 & 1 \\
	\hline
	\multirow{2}{*}{\texttt{bndEdges\{2\}}} & 2 & 3 \\
	                                        & 3 & 4 \\
	\hline
	\texttt{boundaries\{1\}}                & 1 & 2 \\
	\hline
	\texttt{boundaries\{2\}}                & 4 & 6 \\
\end{tabular}}
	\caption{Example mesh on the unit square $(0,1)^2 \subset \R^2$ as well as corresponding data structures.
	Boundary part $1$ (e.g., $\Gamma_D$) is marked in red, boundary part $2$ (e.g., $\Gamma_N$) is marked in green.}
	\label{fig:mesh-example}
\end{figure}
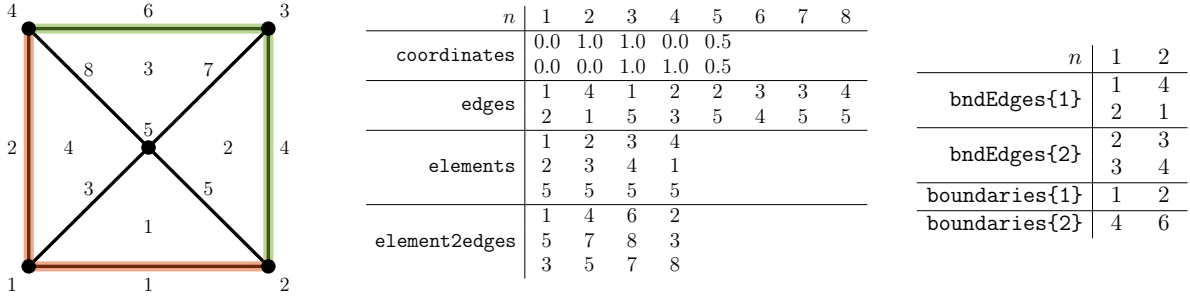

Mesh construction and data structures are showcased in Figure~\ref{fig:mesh-example}.
Note that the orientation of the user-provided edges in \lstinline|bndEdges| is preserved by the automatically generated field \lstinline|edges|.

\subsection{Array layout}\label{subsec:memory-layout}
\begin{figure}
	\centering
	\begin{tikzpicture}[
	mymatrix/.style={
		matrix of nodes,
  		draw,very thick,
  		inner sep=0pt,
  		column sep=-0.2pt,row sep=-0.2pt,
  		nodes in empty cells,
  		cells={
  			nodes={minimum width=1.9em,minimum height=1.9em,draw,very thin,anchor=center,fill=white,text=gray}
  		}
  	}
  ]
 \newcommand{\mycolindex}[1]{\ifnum#1=5 N\else #1\fi}
 \newcommand{\myrowindex}[1]{\ifnum#1=5 T\else #1\fi}
 \matrix[mymatrix,xshift=3em,yshift=3em](matz){
 49 & 52 & 55 & 58 & 61 & 64 & 67 & 70 \\
 50 & 53 & 56 & 59 & 62 & 65 & 68 & 71 \\
 51 & 54 & 57 & 60 & 63 & 66 & 69 & 72 \\
 };
 \matrix[mymatrix,xshift=1.5em,yshift=1.5em](maty){
 25 & 28 & 31 & 34 & 37 & 40 & 43 & 46 \\
 26 & 29 & 32 & 35 & 38 & 41 & 44 & 47 \\
 27 & 30 & 33 & 36 & 39 & 42 & 45 & 48 \\
 };
 \matrix[mymatrix](matx){
 1 & 4 & 7 & 10 & 13 & 16 & 19 & 22 \\
 2 & 5 & 8 & 11 & 14 & 17 & 20 & 23 \\
 3 & 6 & 9 & 12 & 15 & 18 & 21 & 24 \\
 };
 \draw[thick,-stealth] ([xshift=1ex]matx.south east) -- ([xshift=1ex]matz.south east)
  node[text width=5em,xshift=3em,midway,below] {barycentric coordinates};
 \draw[thick,-stealth] ([yshift=-1ex]matx.south west) -- 
  ([yshift=-1ex]matx.south east) node[midway,below] {FEM units};
 \draw[thick,-stealth] ([xshift=-1ex]matx.north west)
   -- ([xshift=-1ex]matx.south west) node[midway,above,rotate=90] {vector data};
 \draw[thick,lightgray] (matx.north west) -- (matz.north west);
 \draw[thick,lightgray] (matx.north east) -- (matz.north east);
 \draw[thick,lightgray] (matx.south east) -- (matz.south east);
\end{tikzpicture}
	\caption{Illustration of the memory layout chosen in our implementation.
	The numbers indicate the order in which the items are stored in memory.}
	\label{fig:memory-layout}
\end{figure}
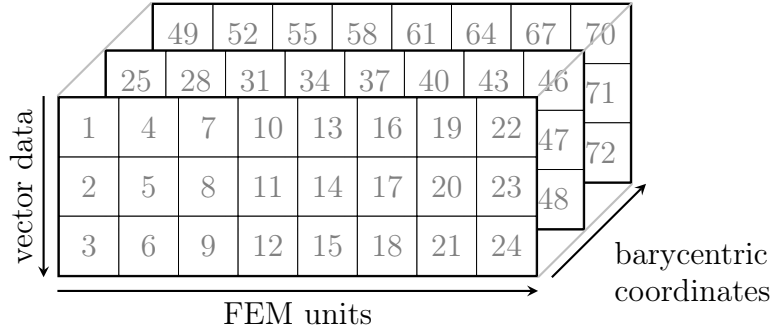

The array layout is chosen such that the first three dimensions of arrays correspond to modeling concepts of finite elements; see Figure~\ref{fig:memory-layout}:
\begin{itemize}
	\item \textbf{1.~dimension (columns)}: This corresponds to the components of vector- or matrix-valued data.
	Matrices are stored column-wise.
	
	\item \textbf{2.~dimension (rows)}: This corresponds to the different units of the mesh, i.e., elements, edges, or vertices.
	
	\item \textbf{3.~dimension (pages)}: This corresponds to different barycentric coordinates, e.g., for numerical quadrature.
\end{itemize}
The rationale behind this order is twofold.
First, objects on the same element often need to be multiplied together.
Hence, it is of advantage if the data representing these objects are continuous in memory.
This is achieved by arranging them along the columns of the array, since \matlab stores arrays in column-major order.
Second, arrays that extend into the third dimension are somewhat clumsy to work with and hard to debug for programmers.
Since evaluations on multiple barycentric coordinates occur mostly internally (e.g., for quadrature rule or plotting), arranging different barycentric coordinates along the third dimension minimizes exposure of the user to more-than-two dimensional arrays.
Within \moofem, one can use the enumeration class \lstinline|Dim| to access these dimensions by \lstinline|Dim.Vector|, \lstinline|Dim.Elements|, and \lstinline|Dim.QuadratureNodes|, respectively.

As an illustrative example, consider the call
\begin{lstlisting}
f = MeshFunction(mesh, @(x) x);
val = eval(f, bary);
\end{lstlisting}
which evaluates $f \colon \Omega \to \R^2 \colon x \mapsto x$ on a collection \lstinline|bary| of barycentric coordinates and a given mesh element-wise.
The value stored in \lstinline|val(i,j,k)| corresponds to $x_i(\lambda^{(k)}, T_j)$, i.e., the $i$-th component of \lstinline|f| evaluated at the $k$-th barycentric coordinate on the $j$-th element.
The matrix valued function
\begin{equation*}
	f \colon \Omega \to \R^{2 \times 2} \colon x \mapsto
	\left(\begin{matrix}
		1 x_1 & 3 x_1 \\
		2 x_1 & 4 x_1
	\end{matrix}\right)
\end{equation*}
can be implemented by \lstinline|f = MeshFunction(mesh, @(x) [1;2;3;4].*x(1,:,:))|.
The corresponding evaluation \lstinline|val(i,j,k)| is equal to $i \cdot x_1(\lambda^{(k)}, T_j)$, since matrices are stored column-wise.
See also Figure~\ref{fig:memory-layout} for a sketch of this memory layout.

\subsection{Efficient linear algebra}\label{subsec:linalg-core}

According to the last section, 3D arrays are essentially interpreted as collections of matrices stored column-wise.
To efficiently execute matrix operations within this memory layout, the function
\begin{lstlisting}
C = vectorProduct(A, B, sizeA, sizeB);
\end{lstlisting}
is used.
It computes the product of two 3D arrays \lstinline|A| and \lstinline|B|, where the first dimension is  interpreted as matrix with given size \lstinline|sizeA| and \lstinline|sizeB|, respectively.
In particular, for all admissible $i,j \in \N$, the output of above call satisfies
\begin{lstlisting}
C(:,i,j) = reshape(A(:,i,j), sizeA) * reshape(B(:,i,j), sizeB);	
\end{lstlisting}
If either \lstinline|sizeA| or \lstinline|sizeB| is a column-vector, the corresponding factor in the above listing is transposed.
For an example, consider the vector $\texttt{v} := (1,2,3,4,5,6)^\intercal$, which, in our memory model, can be interpreted as 
\begin{equation*}
	\texttt{[2,3]}: \quad
	A := \left(\begin{matrix}
		1 & 3 & 5 \\
		2 & 4 & 6
	\end{matrix}\right)
	\quad \text{ or } \quad
	\texttt{[3,2]}: \quad
	B := \left(\begin{matrix}
	1 & 4 \\
	2 & 5 \\
	3 & 6
	\end{matrix}\right).
\end{equation*}
Clearly there holds $A^\intercal \neq B$.
Therefore, transposing the size is necessary to indicate transposition of the corresponding factor in the matrix product:
\begin{align*}
	A B &= \texttt{vectorProduct(v, v, [2,3], [3,2])},\\
	A A^\intercal &= \texttt{vectorProduct(v, v, [2,3], [2,3]')},\\
	B^\intercal B &= \texttt{vectorProduct(v, v, [3,2]', [3,2])}.
\end{align*}

The function \lstinline|vectorProduct| thus enables all possible matrix operations within our memory layout in a convenient, yet very efficient way.
In fact, this routine can also deal with arrays of arbitrary dimension, where extension of singleton dimensions is done automatically, as is common in \matlab.
Per default, the sizes are chosen such that the dot product
\begin{lstlisting}
C(:,i,j) = A(:,i,j)' * B(:,i,j);	
\end{lstlisting}
is computed, which is the most common application; this is also reflected in the name.


\section{Examples}\label{sec:examples}

In this section, we discuss several extensions of the basic AFEM algorithm that is implemented in Listing~\ref{lst:afem}.
We do not claim that \moofem can deal with all FEM applications out of the box, but are convinced that our code structure makes extensions and modifications relatively easy.

\subsection{Higher order AFEM with known solution}\label{subsec:ho-afem}

As a first example we consider the L-shape $\Omega := (-1,1)^2 \backslash \big([0,1] \times [-1,0]\big)$ with boundary parts $\Gamma_R = \emptyset$, $\Gamma_N := \big([0,1] \times \{0\}\big) \cup \big(\{0\} \times [-1,0]\big)$, and $\Gamma_D := \partial \Omega \backslash \Gamma_N$.
With $\big(r(x), \varphi(x)\big)$ being the polar coordinates of $x \in \R^2$, we prescribe the exact solution
\begin{equation*}
\label{eq:uex}
	u(x)
	:=
	r(x)^{2/3} \sin(2\pi/3)
\end{equation*}
and note that it solves
\begin{equation}
\label{eq:example1-model}
	-\Delta u = 0 ~\text{ in } \Omega, \quad
	\nabla u \cdot \normalvec =: \phi ~\text{ on } \Gamma_N, \quad
	u = 0 ~\text{ on } \Gamma_D.
\end{equation}

To solve~\eqref{eq:example1-model} numerically with \moofem, only minor adjustments of the code from Listing~\ref{lst:afem} are necessary.
First, obviously, the correct mesh must be loaded via
\begin{lstlisting}
mesh = Mesh.loadFromGeometry('Lshape');
\end{lstlisting}
This geometry has two predefined boundary parts: The first (boundary index 1) is at the re-entrant corner ($\Gamma_D$), the second (boundary index 2) is everything else ($\Gamma_N$).
Next, the finite element space has to be chosen accordingly with some $p \in \N$:
\begin{lstlisting}
fes = FeSpace(mesh, HigherOrderH1Fe(p), 'dirichlet', 1);
\end{lstlisting}
This loads an implementation of $\SS^p_D(\TT_0)$.
No further adjustments regarding the implementation of higher-order finite elements are necessary.

The next changes concern the coefficients of the linear form \lstinline|lf|.
Our problem~\eqref{eq:example1-model} includes only the Neumann part of the right-hand side from~\eqref{eq:discrete-system}, which can be implemented by the following listing.
\begin{lstlisting}
lf.neumann = MeshFunction(mesh, @exactSolutionNeumannData);
lf.bndNeumann = 2;
function y = exactSolutionNeumannData(x)
	x1 = x(1,:,:);
	x2 = x(2,:,:);
	% determine boundary parts
	right  = (x1 > 0) & (abs(x1) > abs(x2));
	left   = (x1 < 0) & (abs(x1) > abs(x2));
	top    = (x2 > 0) & (abs(x1) < abs(x2));
	bottom = (x2 < 0) & (abs(x1) < abs(x2));
	% compute d/dn u
	[phi, r] = cart2pol(x(1,:,:), x(2,:,:));
	Cr = 2/3 * r.^(-4/3);
	Cphi = 2/3 * (phi + 2*pi*(phi < 0));
	dudx = Cr .* (x1.*sin(Cphi) - x2.*cos(Cphi));
	dudy = Cr .* (x2.*sin(Cphi) + x1.*cos(Cphi));
	y = zeros(size(x1));
	y(right)  = dudx(right);
	y(left)   =-dudx(left);
	y(top)    = dudy(top);
	y(bottom) =-dudy(bottom);
end
\end{lstlisting}
Here, the main part is the implementation of the Neumann derivative $\nabla u \cdot \normalvec$, rather than making the function available to the assembly routines.
Finally, we need to set quadrature rules of sufficiently high order for the corresponding terms of the (bi-)linear form:
\begin{lstlisting}
blf.qra = QuadratureRule.ofOrder(max(2*p-2, 1));
lf.qrNeumann = QuadratureRule.ofOrder(2*p, '1D');
\end{lstlisting}
With these preparatory steps, the FEM system is solved by lines 10--13 of Listing~\ref{lst:afem}.

Finally, the \textsl{a posteriori} indicators have to be adjusted to the current setting:
\begin{equation}
\label{eq:example1-aposteriori}
	\eta_\coarse(T)^2
	=
	h_T^2 \norm{\Delta u_\coarse}_{L^2(T)}^2
	+ h_T \big[ \norm{\jump{\nabla u_\coarse \cdot \normalvec}}_{L^2(\partial T \cap \Omega)}^2
		+ \norm{(\nabla u_\coarse - \nabla u) \cdot \normalvec}_{L^2(\partial T \cap \Gamma_N)}^2 \big].
\end{equation}
Recall from Section~\ref{subsec:memory-layout} that matrices are stored column-wise in the first dimension of 3D arrays.
Thus, the $L^2$-norm of the volume term can be computed by
\begin{lstlisting}
f = CompositeFunction(@(D2u) (D2u(1,:,:) + D2u(4,:,:)).^2, Hessian(u));
qr = QuadratureRule.ofOrder(max(2*(p-2), 1));
volumeRes = integrateElement(f, qr);
\end{lstlisting}
The edgewise $L^2$-norms are a bit more involved.
This is handled by
\begin{lstlisting}
qr = QuadratureRule.ofOrder(p, '1D');
edgeRes = integrateNormalJump(Gradient(u), qr, ...
	@(j) zeros(size(j)), {}, mesh.boundaries{1}, ...
	@(j,phi) j-phi, {lf.neumann}, mesh.boundaries{2}, ...
	@(j) j.^2, {}, ':');
\end{lstlisting}
The syntax of \lstinline|integrateNormalJump| is explained in Section~\ref{subsec:module-integration}.
First, the jump $\jump{\nabla u_\coarse \cdot \normalvec}$ is computed on every edge.
Then, since the edges on $\Gamma_D$ (boundary index 1) do not contribute to the error estimator, the corresponding contributions are set to zero.
Furthermore, the term $\nabla u \cdot \normalvec$, which is stored in \lstinline|lf.neumann|, is subtracted on $\Gamma_N$ (boundary index 2).
Finally, every edge contribution is squared to obtain the edgewise $L^2$-norms of~\eqref{eq:example1-aposteriori}.

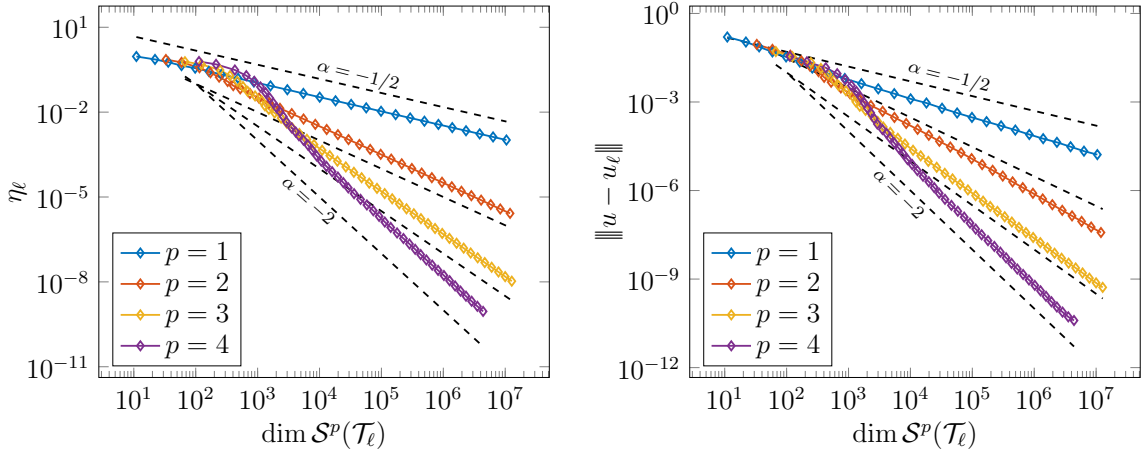
\begin{figure}
	\centering
	\resizebox{0.95\textwidth}{!}{\tikzstyle{plot}=[thick,mark=diamond]%
\tikzstyle{reference}=[thick,dashed]%
\pgfplotstableread[col sep = comma]{figures/data/afemP1.dat}{\afemPone}%
\pgfplotstableread[col sep = comma]{figures/data/afemP2.dat}{\afemPtwo}%
\pgfplotstableread[col sep = comma]{figures/data/afemP3.dat}{\afemPthree}%
\pgfplotstableread[col sep = comma]{figures/data/afemP4.dat}{\afemPfour}%

\begin{tikzpicture}
	
\begin{loglogaxis}[xlabel={$\dim \SS^p(\TT_\ell)$}, ylabel={$\eta_\ell$}, legend pos=south west]

	\addplot[reference,forget plot] table [x={nDofs}, y expr={1.5e1/sqrt(\thisrowno{1})}] {\afemPone} node[below,rotate=-14] at (axis cs:5E4,9E-1) {\tiny{$\alpha = -1/2$}};
	\addplot[reference,forget plot] table [x={nDofs}, y expr={1e1/\thisrowno{1}}] {\afemPthree};
	\addplot[reference,forget plot] table [x={nDofs}, y expr={1e2/(sqrt(\thisrowno{1})*\thisrowno{1})}] {\afemPthree};
	\addplot[reference,forget plot] table [x={nDofs}, y expr={1e3/(\thisrowno{1}*\thisrowno{1})}] {\afemPfour} node[below,rotate=-42] at (axis cs:1E4,2E-5) {\tiny{$\alpha = -2$}};

	\addplot[plot,color1] table [x={nDofs}, y={estimator}] {\afemPone}; \addlegendentry{$p=1$}
	\addplot[plot,color2] table [x={nDofs}, y={estimator}] {\afemPtwo}; \addlegendentry{$p=2$}
	\addplot[plot,color3] table [x={nDofs}, y={estimator}] {\afemPthree}; \addlegendentry{$p=3$}
	\addplot[plot,color4] table [x={nDofs}, y={estimator}] {\afemPfour}; \addlegendentry{$p=4$}
\end{loglogaxis}

\end{tikzpicture}
\quad
\begin{tikzpicture}

\begin{loglogaxis}[xlabel={$\dim \SS^p(\TT_\ell)$}, ylabel={$\enorm{u-u_\ell}$}, legend pos=south west]

	\addplot[reference,forget plot] table [x={nDofs}, y expr={5e-1/sqrt(\thisrowno{1})}] {\afemPone} node[below,rotate=-14] at (axis cs:5E4,3E-2) {\tiny{$\alpha = -1/2$}};
	\addplot[reference,forget plot] table [x={nDofs}, y expr={3e0/\thisrowno{1}}] {\afemPthree};
	\addplot[reference,forget plot] table [x={nDofs}, y expr={1e1/(sqrt(\thisrowno{1})*\thisrowno{1})}] {\afemPthree};
	\addplot[reference,forget plot] table [x={nDofs}, y expr={1e2/(\thisrowno{1}*\thisrowno{1})}] {\afemPfour} node[below,rotate=-45] at (axis cs:1E4,2E-6) {\tiny{$\alpha = -2$}};
	
	\addplot[plot,color1] table [x={nDofs}, y={H1Error}] {\afemPone}; \addlegendentry{$p=1$}
	\addplot[plot,color2] table [x={nDofs}, y={H1Error}] {\afemPtwo}; \addlegendentry{$p=2$}
	\addplot[plot,color3] table [x={nDofs}, y={H1Error}] {\afemPthree}; \addlegendentry{$p=3$}
	\addplot[plot,color4] table [x={nDofs}, y={H1Error}] {\afemPfour}; \addlegendentry{$p=4$}
\end{loglogaxis}

\end{tikzpicture}}
	\caption{Error estimator $\eta_\ell$ (left) and energy error $\enorm{u-u_\ell}$ (right) over number of DOFs for problem~\eqref{eq:example1-model} with different polynomial orders $p$.}
	\label{fig:example1}
\end{figure}

\begin{figure}
	\centering
	\resizebox{0.95\textwidth}{!}{\tikzstyle{plot}[diamond]=[thick,mark=#1]%
\tikzstyle{reference}=[thick,dashed]%
\pgfplotstableread[col sep = comma]{figures/data/timingP1.dat}{\timingPone}%
\pgfplotstableread[col sep = comma]{figures/data/timingP4.dat}{\timingPfour}%

\begin{tikzpicture}
	
\begin{loglogaxis}[xlabel={$\dim \SS^p(\TT_\ell)$}, ylabel={computation time [s] ${} / \dim \SS^p(\TT_\ell)$},%
	ymax=2e-1,grid=both,grid style={line width=.1pt, draw=gray!10},major grid style={line width=.2pt,draw=gray!50},%
	legend style={legend pos=north east, anchor=north east, font=\tiny, legend columns=2}]

	\addplot[plot,color1] table [x={nDofs}, y expr={\thisrowno{4}/\thisrowno{1}}] {\timingPone}; \addlegendentry{solve}
	\addplot[plot=x,color1] table [x={nDofs}, y expr={\thisrowno{5}/\thisrowno{1}}] {\timingPone}; \addlegendentry{estimate}
	\addplot[plot,color2] table [x={nDofs}, y expr={\thisrowno{2}/\thisrowno{1}}] {\timingPone}; \addlegendentry{assemble A}
	\addplot[plot=x,color2] table [x={nDofs}, y expr={\thisrowno{6}/\thisrowno{1}}] {\timingPone}; \addlegendentry{mark}
	\addplot[plot,color3] table [x={nDofs}, y expr={\thisrowno{7}/\thisrowno{1}}] {\timingPone}; \addlegendentry{refine}
	\addplot[plot=x,color3] table [x={nDofs}, y expr={\thisrowno{3}/\thisrowno{1}}] {\timingPone}; \addlegendentry{assemble F}
	\addplot[plot=o,color4] table [x={nDofs}, y expr={\thisrowno{8}/\thisrowno{1}}] {\timingPone}; \addlegendentry{all combined}
\end{loglogaxis}

\end{tikzpicture}
\quad
\begin{tikzpicture}

\begin{loglogaxis}[xlabel={$\dim \SS^p(\TT_\ell)$},%
	ymin=2e-9,ymax=9e-2,grid=both,grid style={line width=.1pt, draw=gray!10},major grid style={line width=.2pt,draw=gray!50},%
	legend style={legend pos=north east, anchor=north east, font=\tiny, legend columns=2}]

	\addplot[plot,color1] table [x={nDofs}, y expr={\thisrowno{4}/\thisrowno{1}}] {\timingPfour}; \addlegendentry{solve}
	\addplot[plot=x,color1] table [x={nDofs}, y expr={\thisrowno{5}/\thisrowno{1}}] {\timingPfour}; \addlegendentry{estimate}
	\addplot[plot,color2] table [x={nDofs}, y expr={\thisrowno{2}/\thisrowno{1}}] {\timingPfour}; \addlegendentry{assemble A}
	\addplot[plot=x,color2] table [x={nDofs}, y expr={\thisrowno{6}/\thisrowno{1}}] {\timingPfour}; \addlegendentry{mark}
	\addplot[plot,color3] table [x={nDofs}, y expr={\thisrowno{7}/\thisrowno{1}}] {\timingPfour}; \addlegendentry{refine}
	\addplot[plot=x,color3] table [x={nDofs}, y expr={\thisrowno{3}/\thisrowno{1}}] {\timingPfour}; \addlegendentry{assemble F}
	\addplot[plot=o,color4] table [x={nDofs}, y expr={\thisrowno{8}/\thisrowno{1}}] {\timingPfour}; \addlegendentry{all combined}
\end{loglogaxis}

\end{tikzpicture}}
	\caption{Computation time per DOF over number of DOFs for the different parts of the AFEM algorithm for problem~\eqref{eq:example1-model} with polynomial degree $p=1$ (left) and $p=4$ (right).}
	\label{fig:timing}
\end{figure}
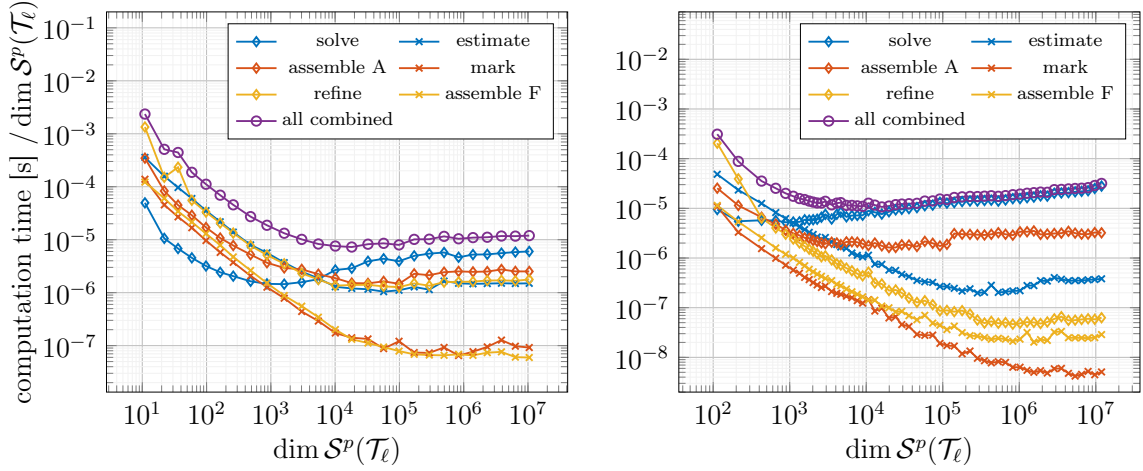

The remainder of the code is analogous to the one presented in Listing~\ref{lst:afem}.
Note that virtually all changes are due to the different model problem and not for implementational reasons.
Figure~\ref{fig:example1} shows the results obtained for $p=1,2,3,4$.
Note that, from an implementational point of view, the polynomial degree $p \in \N$ can be chosen arbitrarily high.
Computation times for the different parts of the adaptive algorithm are shown in Figure~\ref{fig:timing}.
In both the lowest and the higher order case, most time is spent for solution of the linear system.
In the higher order case, one clearly sees that solving with the \matlab backslash operator has more than linear complexity.

As a final note, the exact error of the finite element solution $u_\coarse$ can be easily computed by the following code snippet.
Recall that \lstinline|A| is the finite element matrix of the Laplacian.
\begin{lstlisting}
uex = FeFunction(fes);
uex.setData(nodalInterpolation(MeshFunction(mesh, @exactSolution), fes));
deltaU = u.data - uex.data;
H1Error = sqrt(deltaU * A * deltaU');

function y = exactSolution(x)
	[phi, r] = cart2pol(x(1,:,:), x(2,:,:));
	phi = phi + 2*pi*(phi < 0);
	y = r.^(2/3) .* sin(2/3 * phi);
end
\end{lstlisting}

\subsection{Goal-oriented AFEM with discontinuous data}\label{subsec:goafem}

With $\Omega := (0,1)^2$ and $\Gamma_D := \partial \Omega$, we consider an example from~\cite{ms2009}:
\begin{equation}
\label{eq:example2-model}
	-\Delta u = -\div\vec{f} ~\text{ in } \Omega, \quad
	u = 0 ~\text{ on } \Gamma_D, \quad
	\text{where }
	\vec{f}(x)
	:=
	\begin{cases}
		(1,0) & \text{ if } x_1 + x_2 < 1/2,\\
		(0,0) & \text{ else.}
	\end{cases}
\end{equation}
For most FEM software, discontinuous coefficients or data demand some caution: for quadrature nodes that lie on the discontinuity, evaluation is not well-defined.
A first solution is to make the initial triangulation $\TT_0$ of $\Omega$ resolve the regions of discontinuity.
In our case, this can be achieved by uniform refinement using the RGB-strategy:
\begin{lstlisting}
mesh = Mesh.loadFromGeometry('unitsquare');
mesh.refineUniform(1, 'RGB');
\end{lstlisting}
This is also needed for residual error estimators, since they are comprised of element-wise $L^2$-norms of $\div\vec{f}$ (which vanishes if the discontinuity is resolved by the mesh and is not defined otherwise).

A second problem is the jump term $\jump{\cdot}$ in the error estimators, since this is evaluated on edges, where the discontinuity now lies.
This can be solved by interpolating the data to a non-continuous FEM space.
To obtain vector-valued data, we first interpolate the non-continuous first component and then compose this with the vanishing second component, according to the memory layout presented in Section~\ref{subsec:memory-layout}:
\begin{lstlisting}
ncFes = FeSpace(mesh, LowestOrderL2Fe);
w = FeFunction(ncFes);
chiT = MeshFunction(mesh, @(x) sum(x, Dim.Vector) < 1/2);
w.setData(nodalInterpolation(chiT, ncFes));
lfF = LinearForm(fes);
lfF.fvec = CompositeFunction(@(w) [w;zeros(size(w))], w);
\end{lstlisting}
The nodal interpolation in the listing above only sets the data for \lstinline|w| on the initial mesh $\TT_0$.
To have \lstinline|w| available on refined meshes, we can repeat this interpolation process after every mesh refinement.
A more efficient method is to use the prolongation class \lstinline|P = LoFeProlongation(fes)| that is tailored specifically to lowest order $L^2$- and $H^1$-elements; see Section~\ref{subsec:module-fem}.
Data for prolongation is computed automatically whenever the mesh is refined; see Section~\ref{subsec:module-geometry}.
After updating \lstinline|w| by \lstinline|w.setData(P.prolongate(w))|, the next call of \lstinline|assemble(lfF)| already yields the updated right-hand side, since the coefficient \lstinline|lfF.fvec| stores a reference to \lstinline|w|.

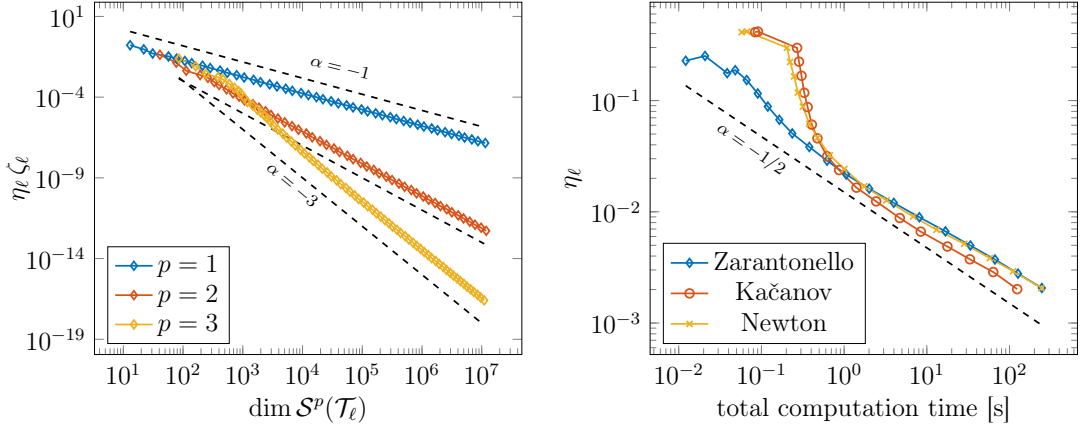
\begin{figure}
	\centering
	\resizebox{0.9\textwidth}{!}{\begin{tikzpicture}
\tikzstyle{plot}=[thick,mark=diamond]%
\tikzstyle{reference}=[thick,dashed]%
\pgfplotstableread[col sep = comma]{figures/data/goafemP1.dat}{\goafemPone}%
\pgfplotstableread[col sep = comma]{figures/data/goafemP2.dat}{\goafemPtwo}%
\pgfplotstableread[col sep = comma]{figures/data/goafemP3.dat}{\goafemPthree}%
	
\begin{loglogaxis}[xlabel={$\dim \SS^p(\TT_\ell)$}, ylabel={$\eta_\ell \, \zeta_\ell$}, legend pos=south west]

	\addplot[reference,forget plot] table [x={nDofs}, y expr={1.5e1/\thisrowno{1}}] {\goafemPone} node[below,rotate=-14] at (axis cs:5E4,6E-2) {\tiny{$\alpha = -1$}};
	\addplot[reference,forget plot] table [x={nDofs}, y expr={1e1/(\thisrowno{1}*\thisrowno{1})}] {\goafemPthree};
	\addplot[reference,forget plot] table [x={nDofs}, y expr={1e3/(\thisrowno{1}*\thisrowno{1}*\thisrowno{1})}] {\goafemPthree} node[below,rotate=-40] at (axis cs:1E4,2E-9) {\tiny{$\alpha = -3$}};

	\addplot[plot,color1] table [x={nDofs}, y={goalErrorEstimate}] {\goafemPone}; \addlegendentry{$p=1$}
	\addplot[plot,color2] table [x={nDofs}, y={goalErrorEstimate}] {\goafemPtwo}; \addlegendentry{$p=2$}
	\addplot[plot,color3] table [x={nDofs}, y={goalErrorEstimate}] {\goafemPthree}; \addlegendentry{$p=3$}
\end{loglogaxis}

\end{tikzpicture}\quad\tikzstyle{plot}=[thick]%
\tikzstyle{reference}=[thick,dashed]%
\pgfplotstableread[col sep = comma]{figures/data/ailfem-zarantonello.dat}{\zarantonello}%
\pgfplotstableread[col sep = comma]{figures/data/ailfem-kacanov.dat}{\kacanov}%
\pgfplotstableread[col sep = comma]{figures/data/ailfem-newton.dat}{\newton}%

\begin{tikzpicture}
	
\begin{loglogaxis}[xlabel={total computation time [s]}, ylabel={$\eta_\ell$}, legend pos=south west]

	\addplot[reference,forget plot] table [x={time}, y expr={1.5e-2/sqrt(\thisrowno{3})}] {\zarantonello} node[below,rotate=-35] at (axis cs:1E-1,5E-2) {\tiny{$\alpha = -1/2$}};

	\addplot[plot,color1,mark=diamond] table [x={time}, y={estimator}] {\zarantonello}; \addlegendentry{Zarantonello}
	\addplot[plot,color2,mark=o] table [x={time}, y={estimator}] {\kacanov}; \addlegendentry{Ka\v{c}anov}
	\addplot[plot,color3,mark=x] table [x={time}, y={estimator}] {\newton}; \addlegendentry{Newton}
\end{loglogaxis}

\end{tikzpicture}}
	\caption{Left: Estimator for the goal error~\eqref{eq:goal-error-estimate} over number of DOFs for problem~\eqref{eq:example2-model} from Section~\ref{subsec:goafem} with different polynomial orders $p$.
	Right: Error estimator over total computation time for the linearization methods from Section~\ref{subsec:nonlinear}.}
	\label{fig:example2+3}
\end{figure}

\begin{figure}
	\centering
	\includegraphics[width=0.35\textwidth]{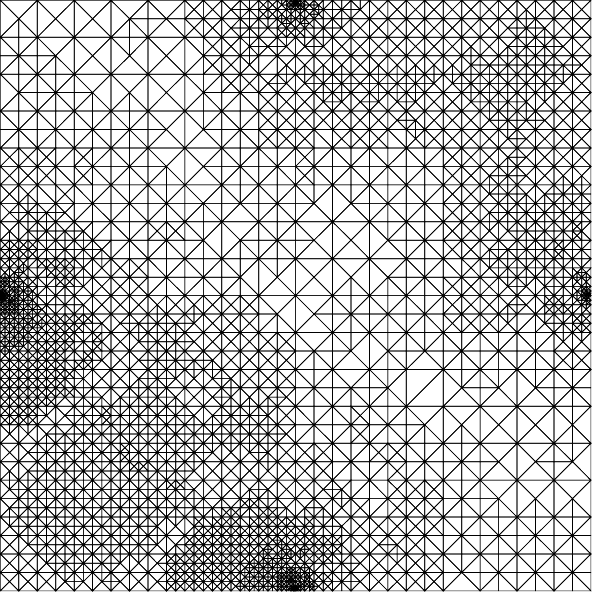}
	\qquad
	\includegraphics[width=0.35\textwidth]{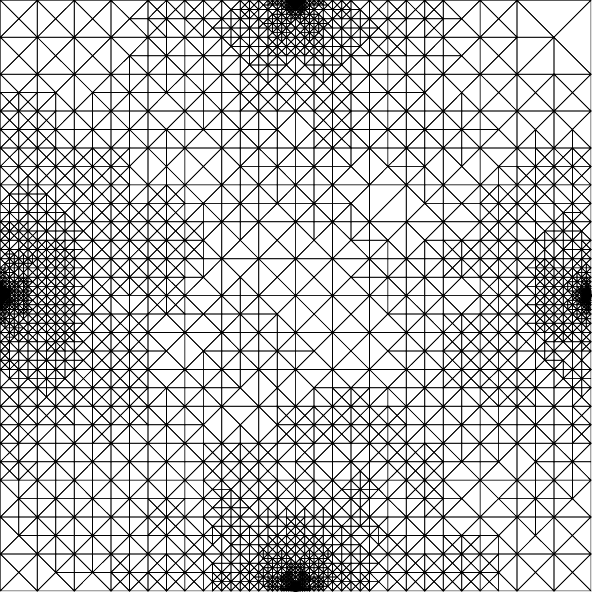}
	\caption{Meshes generated from the GOAFEM algorithm from Section~\ref{subsec:goafem} with polynomial orders $p=1$ (left) and $p=3$ (right).}
	\label{fig:mesh-goafem}
\end{figure}

In \emph{goal-oriented} adaptive FEM (GOAFEM), we are interested in the \emph{goal value} $G(u)$ for a linear functional
\begin{equation*}
	 G \colon
	 H^1_D(\Omega) \to \R,
	 \quad
	 G(v)
	 =
	 \int_\Omega \vec{g} \cdot \nabla v \d{x}
	 \quad \text{ with } \quad
	 \vec{g}(x)
	 :=
	 \begin{cases}
		 (-1,0) & \text{ if } x_1 + x_2 > 3/2,\\
		 (0,0) & \text{ else.}
	 \end{cases}
\end{equation*}
Approximating the goal value $G(u)$ is often more interesting in applications than approximating the solution $u$ as a whole.
To efficiently approximate the goal value in the spirit of Algorithm~\ref{alg:afem}, one introduces the so-called \emph{dual} problem
\begin{equation}
\label{eq:example2-dual}
	-\Delta z = -\div\vec{g} ~\text{ in } \Omega, \quad
	z = 0 ~\text{ on } \Gamma_D,
\end{equation}
which can be implemented and solved analogously to~\eqref{eq:example2-model}.
Since solving is often the most time consuming part of AFEM, we can do this in parallel for~\eqref{eq:example2-model} and~\eqref{eq:example2-dual}:
\begin{lstlisting}
rhs = [assemble(lfF), assemble(lfG)];
uz = A(freeDofs,freeDofs) \ rhs(freeDofs,:);
u.setFreeData(uz(:,1));
z.setFreeData(uz(:,2));
\end{lstlisting}

After solving~\eqref{eq:example2-model} and~\eqref{eq:example2-dual} by FEM on a triangulation $\TT_\coarse$ to obtain the discrete solutions $u_\coarse$ and $z_\coarse$, respectively, one can compute the \textsl{a posteriori} residual error estimators
\begin{align*}
\label{eq:example2-aposteriori}
	\eta_\coarse(T)^2
	&=
	h_T^2 \norm{\Delta u_\coarse}_{L^2(T)}^2
	+ h_T \norm{\jump{(\nabla u_\coarse - \vec{f}) \cdot \normalvec}}_{L^2(\partial T \cap \Omega)}^2,\\
	\zeta_\coarse(T)^2
	&=
	h_T^2 \norm{\Delta z_\coarse}_{L^2(T)}^2
	+ h_T \norm{\jump{(\nabla z_\coarse - \vec{g}) \cdot \normalvec}}_{L^2(\partial T \cap \Omega)}^2
\end{align*}
analogously to~\eqref{eq:example1-aposteriori}.
The error in the goal functional is controlled by the estimator product
\begin{equation}
\label{eq:goal-error-estimate}
	|G(u) - G(u_\coarse)|
	\lesssim
	\Big[ \sum_{T \in \TT_\coarse} \eta_\coarse(T)^2 \Big]^{1/2}
	\Big[ \sum_{T \in \TT_\coarse} \zeta_\coarse(T)^2 \Big]^{1/2},
\end{equation}
for which different marking criteria have been analyzed~\cite{fpz2016}.
Thus, the remaining implementation comprises only minor modifications of Listing~\ref{lst:afem}.
The upper bounds of this last equation for different polynomial orders $p$ can be seen in Figure~\ref{fig:example2+3} and the resulting meshes for $p=1,3$ are shown in Figure~\ref{fig:mesh-goafem}.

\subsection{Iterative solution of nonlinear equations}\label{subsec:nonlinear}

In this last example, we consider the L-shape $\Omega := (-1,1)^2 \backslash \big([0,1] \times [-1,0]\big)$ with Dirichlet boundary $\Gamma_D := \partial \Omega$.
On this domain, we consider the quasi-linear problem
\begin{equation}
\label{eq:example3-model}
	-\div \Big( \mu \big(|\nabla u|^2 \big) \nabla u \Big)
	= 1 ~\text{ in } \Omega, \quad
	u = 0 ~\text{ on } \Gamma_D,
	\quad \text{with} \quad
	\mu(t) = 1 + \exp(-t).
\end{equation}
This is a variation of an example given in~\cite{hpw2021}, where also adaptive iterative linearization techniques (AILFEM) for this class of problems are presented.
With a given initial guess $u^0 \in H^1_D(\Omega)$, we consider the following linearizations:
\begin{enumerate}[label={(\roman*)}]
	\item \textbf{Zarantonello iteration}:
	Let $\delta > 0$ be sufficiently small.
	Given $u^n \in H^1_D(\Omega)$, the next iterate $u^{n+1} \in H^1_D(\Omega)$ reads $u^{n+1} := u^n + \delta v$, where $v \in H^1_D(\Omega)$ solves
	\begin{equation*}
	\label{eq:zarantonello}
		-\Delta v
		=
		\div \Big( \mu \big(|\nabla u^n|^2 \big) \nabla u^n \Big) + 1.
	\end{equation*}
	
	\item \textbf{Ka\v{c}anov iteration}:
	Given $u^n \in H^1_D(\Omega)$, the next iterate $u^{n+1} \in H^1_D(\Omega)$ solves
	\begin{equation*}
	\label{eq:kacanov}
		-\div \Big( \mu \big(|\nabla u^n|^2 \big) \nabla u^{n+1} \Big) = 1.
	\end{equation*}
	
	\item \textbf{Newton iteration}:
	Given $u^n \in H^1_D(\Omega)$, the next iterate $u^{n+1} \in H^1_D(\Omega)$ reads $u^{n+1} := u^{n} + v$, where $v \in H^1_D(\Omega)$ solves
	\begin{equation*}
	\label{eq:newton}
		-\div \Big( \mu (|\nabla u^n|^2) \nabla v + 2 \mu'(|\nabla u^n|^2) \big( \nabla u^n \otimes \nabla u^n \big) \nabla v \Big)
		=
		\div \Big( \mu \big(|\nabla u^n|^2 \big) \nabla u^n \Big) + 1.
	\end{equation*}
\end{enumerate}
All iterations~(i)--(iii) feature coefficients that depend in a nonlinear fashion on the previous iterate $u^n$.
However, their implementation is relatively simple, owing to the uniform evaluation mechanics of the \lstinline|Evaluable| interface, from which also \lstinline|FeFunction| is derived.
Assuming that, for some triangulation $\TT_\coarse$ of $\Omega$, the previous iterates $u_\coarse^n$ correspond to the \lstinline|FeFunction| instance \lstinline|u|, the following code snippet acts as template for all three iterations with $u_\coarse^0 = 0$:
\begin{lstlisting}
% set coefficients of blf & lf
u = FeFunction(fes);
u.setData(0);
v = FeFunction(fes);
freeDofs = getFreeDofs(fes);
while true
	A = assemble(blf);
	F = assemble(lf);
	% solve linear systems and update data of u
end
\end{lstlisting}
The two steps in this template that are merely outlined in a comment differ for each method.
They are described in the following listing, separated by comments:
\begin{lstlisting}
% --- Zarantonello: setup
blf.a = Constant(mesh, 1);
lf.f = Constant(mesh, 1);
lf.fvec = CompositeFunction(@(p) -mu(vectorProduct(p, p)) .* p, Gradient(u));
% --- Zarantonello: update
v.setFreeData(A(freeDofs,freeDofs) \ F(freeDofs));
u.setData(u.data + delta*v.data);
% --- Kacanov: setup
blf.a = CompositeFunction(@(p) mu(vectorProduct(p, p)), Gradient(u));
lf.f = Constant(mesh, 1);
% --- Kacanov: update
u.setFreeData(A(freeDofs,freeDofs) \ F(freeDofs));
% --- Newton: setup
blf.a = CompositeFunction(@(p) mu(vectorProduct(p, p)) .* [1;0;0;1] ...
	+ 2*muPrime(vectorProduct(p, p)).*vectorProduct(p, p, [2,1], [2,1]'), ...
	Gradient(u));
lf.f = Constant(mesh, 1);
lf.fvec = CompositeFunction(@(p) -mu(vectorProduct(p, p)) .* p, Gradient(u));
% --- Newton: update
v.setFreeData(A(freeDofs,freeDofs) \ F(freeDofs));
u.setData(u.data + v.data);
% --- Additional functions
mu = @(t) 1 + exp(-t);
muPrime = @(t) -exp(-t);
\end{lstlisting}
As explained in Section~\ref{subsec:linalg-core}, with \lstinline|p = Gradient(u)|, the two calls of \lstinline|vectorProduct| in the Newton bilinear form represent
\begin{align*}
	|\nabla u_\coarse^n|^2
	&= \text{\lstinline|vectorProduct(p, p)|},\\
	\nabla u_\coarse^n \otimes \nabla u_\coarse^n
	&= \text{\lstinline|vectorProduct(p, p, [2,1], [2,1]')|}.
\end{align*}

To get an adaptive algorithm in the spirit of Algorithm~\ref{alg:afem} for lowest order FEM, i.e., $p=1$, error estimation is done by
\begin{equation*}
	\eta_\coarse(T)^2
	:=
	h_T^2 \norm{1}_{L^2(T)}^2 + h_T \norm{\jump{\mu(|\nabla u^n_\coarse|^2) \nabla u_\coarse^n \cdot \normalvec}}_{L^2(\partial T \cap \Omega)}^2,
\end{equation*}
which is analogous to~\eqref{eq:example1-aposteriori}.
Finally, we remark that \cite{hpw2021} suggests to use $u_0^0 = 0 \in \SS^1_D(\TT_0)$ only on the coarsest level and then to proceed by nested iteration
\begin{equation*}
	u_{\ell+1}^0
	:=
	u_\ell^{\underline{n}(\ell)} \in \SS^1_D(\TT_{\ell+1})
	\quad
	\text{for all } \ell \in \N,
\end{equation*}
where $\underline{n}(\ell)$ is the last iteration on the previous level $\TT_\ell$, i.e., $u^{\underline{n}(\ell)}_\ell$ is the final iterate on $\TT_\ell$.
For lowest-order $H^1_D(\Omega)$-conforming FEM, this can be done by the prolongation class \lstinline|LoFeProlongation|; see Section~\ref{subsec:module-fem}.

A numerical comparison of the three presented iterative linearization methods can be seen in Figure~\ref{fig:example2+3}.

\printbibliography

\end{document}